\documentclass{article}

\usepackage{amsmath}
\usepackage{amssymb}
\usepackage{graphicx}
\usepackage{bm}
\usepackage{mathdots}
\usepackage{booktabs}
\usepackage{rotating}
\usepackage[ruled,linesnumbered]{algorithm2e}
\usepackage[a4paper,left=2.8cm,right=2.8cm,top=2.5cm,bottom=2.5cm]{geometry}
\usepackage{fancyhdr}
\usepackage{hyperref}
\usepackage{url}
\usepackage{mathtools}
\usepackage{multirow}

\usepackage[framemethod=tikz]{mdframed}

\newtheorem{theorem}{Theorem}[section]

\newtheorem{proposition}{Proposition}[section]
\newtheorem{lemma}{Lemma}[section]
\newtheorem{definition}{Definition}[section]

\newtheorem{remark}{Remark}[section]

\makeatletter % `@' now normal "letter"
\numberwithin{equation}{section}
\makeatother  % `@' is restored as "non-letter"

\newenvironment{proof}{{\noindent\it Proof.}\quad}{\hfill $\square$\\}

\allowdisplaybreaks

\newcommand{\hardhyper}{\mathcal{H}_{n}^{\lambda}}
\newcommand{\lassohyper}{\mathcal{L}_{n}^{\lambda}}

\newcommand{\hyper}{\mathcal{L}_n}
\newcommand{\discreteinner}[1]{\left\langle#1\right\rangle_N}
\newcommand{\inner}[1]{\left\langle#1\right\rangle}
\newcommand{\brackets}[1]{\left(#1\right)}
\newcommand{\squarebrackets}[1]{\left[#1\right]}

\newcommand{\absolute}[1]{\left|#1\right|}

\newcommand{\hardoperator}[1]{\eta_{H}\left(\discreteinner{#1},\lambda\right)}
\newcommand{\hardoperatornew}[1]{\eta_{H}\left(#1,\lambda\right)}

\newcommand{\norminfty}[1]{\|#1\|_{\infty}}
\newcommand{\normtwo}[1]{\|#1\|_2}

\newcommand{\normzero}[1]{\|#1\|_0}

\newcommand{\summ}[2]{\sum\limits_{#1 =1}^{#2}}

\begin{document}
% Stable image reconstruction via the springback-penalized model
\title{Hard thresholding hyperinterpolation over general regions}

\author{Congpei An\footnotemark[1],\footnotemark[2] \text{ and } Jiashu Ran\footnotemark[2],\footnotemark[3]}

\renewcommand{\thefootnote}{\fnsymbol{footnote}}
\footnotetext[1]{School of Mathematics and Statistics, Guizhou University, Guiyang 550025, Guizhou, China (ancp@gzu.edu.cn),}
\footnotetext[2]{School of Mathematics, Southwestern University of Finance and Economics, Chengdu, China (andbachcp@gmail.com, ranjiashuu@gmail.com),}
\footnotetext[3]{Department of Mathematical and Statistical Sciences, University of Alberta, Edmonton, Alberta T6G 2G1, Canada (jran@ualberta.ca)}

% \date{}
%\shorttitle{SHORT TITLE}
%\shortauthor{F.~FIRSTA AND S.~SECONDA}
% \usepackage{keywords}

\maketitle

\begin{abstract}
This paper proposes a novel variant of hyperinterpolation, called hard thresholding hyperinterpolation. This approximation scheme of degree $n$ leverages a hard thresholding operator to filter all hyperinterpolation coefficients, which approximate the Fourier coefficients of a continuous function by a quadrature rule with algebraic exactness $2n$. We prove that hard thresholding hyperinterpolation is the unique solution to an $\ell_0$-regularized weighted discrete least squares approximation problem. Hard thresholding hyperinterpolation is not only idempotent and commutative with hyperinterpolation, but also adheres to the Pythagorean theorem in terms of the discrete (semi) inner product. By the estimate of the reciprocal of Christoffel function, we present the upper bound of the uniform norm of hard thresholding hyperinterpolation operator. Additionally, hard thresholding hyperinterpolation possesses denoising and basis selection abilities akin to Lasso hyperinterpolation. To judge the $L_2$ errors of both hard thresholding and Lasso hyperinterpolations, we propose a criterion that integrates the regularization parameter with the product of noise coefficients and the signs of hyperinterpolation coefficients. Numerical examples on the sphere, spherical triangle and the cube demonstrate the denoising ability of hard thresholding hyperinterpolation.
\end{abstract}



\section{Introduction}\label{sec:introduction}

The concept of hyperinterpolation was introduced by Sloan in 1995 \cite{sloan1995hyperinterpolation}. The hyperinterpolation operator $\mathcal{L}_n$ is defined by replacing Fourier integrals in the $L_2$ orthogonal projection onto polynomial spaces with a discrete measure based on a positive weight quadrature rule that has algebraic precision $2n$. The result of \cite{sloan1995hyperinterpolation} has sparked numerous investigations into finding suitable quadrature rules for different regions, thereby expanding the potential applications of hyperinterpolation \cite{an2022exactness,an2024bypassing,an2023efficient}.

Hyperinterpolation is a powerful tool in high-dimensional approximation \cite{dai2006hyperinterpolation,LeGia2001uniform,lin2021distributed,reimer2003multivariate,Wade2013hyperinterpolation,Wang2017needlet}. It requires only the availability of a positive weight quadrature rule that exactly integrates polynomial of degree $2n$  \cite{sloan1995hyperinterpolation}. This means the function of interest must be sampled on a carefully selected finite set to meet the requirements of the quadrature formula. In most applications, functions are often given by sampled data. However, the modern era of high-throughput data collection creates data with abundant noise. To recover functions from noisy data, An and Wu developed Lasso hyperinterpolation \cite{an2021lasso}, which uses a soft thresholding operator to process all hyperinterpolation coefficients. Although Lasso hyperinterpolation does not retain the projection property or basis invariance of classical hyperinterpolation, it offers an efficient method for basis selection and denoising, leading to a sparse solution. Lasso hyperinterpolation \cite{an2021lasso} is the solution to an $\ell_1$-regularized weighted discrete least squares problem, which can be regarded as a convex relaxation of an $\ell_0$-regularized problem \cite{Foucart2013compressing}. As is known to all, an $\ell_0$-regularized weighted discrete least squares problem is NP-Hard to be solved \cite{bruckstein2009sparse}. However, we indicate that hard thresholding hyperinterpolation  uniquely solves, in closed form, an $\ell_0$-regularized weighted discrete least squares problem (see Theorem \ref{thm:ExactSolution}) under the same conditions of hyperinterpolation. 

Hard thresholding hyperinterpolation exploits a hard thresholding operator to filter all hyperinterpolation coefficients for a given positive regularization parameter. Specifically, it retains hyperinterpolation coefficients whose absolute values exceed the regularization parameter and sets all the others to zero. Since the hard thresholding operator removes small coefficients completely, it can effectively eliminate noise while retaining larger or more important coefficients that capture the essential features of the test function. We obtain some inescapable algebraic and geometric properties of hard thresholding hyperinterpolation: it is idempotent and commutative with hyperinterpolation, but fails to be symmetric; the composition of hard thresholding hyperinterpolation and hyperinterpolation is Hermitian; it satisfies the Pythagorean theorem \cite{stein2011fourier} but is not basis invariant. 

In the error analysis, we treat properties of hard thresholding hyperinterpolation from three aspects. First, we prove that the operator norm of hard thresholding hyperinterpolation does not exceed that of hyperinterpolation, and we derive the $L_2$ errors when the test function is contaminated by noise. Second, utilizing the reciprocal of Christoffel function \cite{sommariva2014multivariate}, we provide the uniform norms for both hard thresholding hyperinterpolation and the classical hyperinterpolation. Third, we analyze the $L_2$ errors from a practical perspective, specifically demonstrating that hard thresholding hyperinterpolation achieves lower $L_2$ errors compared to Lasso hyperinterpolation (see Theorem \ref{thm:priorparameter}). 

In the sequel, an overview of the fundamental concepts of hyperinterpolation and Lasso hyperinterpolation is provided. Section \ref{sec:HTH} investigates some characterizations of hard thresholding hyperinterpolation, i.e., algebraic and geometric properties. In section \ref{sec:erroranalysis}, we give error analysis from $L_2$ norms, uniform norms and practical viewpoints.  Finally,  we apply hard thresholding hyperinterpolation to the sphere, the spherical triangle and the cube to verify our theories. 

\section{Hard thresholding hyperinterpolation}\label{sec:pre}
Before the discussion of hard thresholding hyperinterpolation, we introduce hyperinterpolation and Lasso hyperinterpolation.
\subsection{Hyperinterpolation}
Let $\Omega$ be a compact and smooth Riemannian manifold in $\mathbb{R}^s$ with smooth or empty boundary and the measure $\text{d}\omega$, which satisfies that
\begin{equation*}
\int_{\Omega}\text{d}\omega=V<\infty.
\end{equation*}
Given the $L_2$ inner product
\begin{equation}\label{equ:inner}
    \langle f, g \rangle = \int_{\Omega} fg {\rm{d} \omega} \qquad \forall f,g \in L_2(\Omega),
\end{equation}
and the induced norm $\|f\|_2:=\langle f,f\rangle^{1/2}$, where $L_2(\Omega)$ denotes the Hilbert space of square-integrable functions on $\Omega$. Let $\{\Phi_{\ell}\}_{\ell=1}^{d_n}$ be an orthonormal basis of $\mathbb{P}_{n}(\Omega)$, the space of polynomials of total-degree at most $n$ restricted to $\Omega$, with respect to the measure ${\rm{d}}\omega$ on $\Omega$,  where $ d_n := {\rm{dim}}(\mathbb{P}_{n}(\Omega))$.

Considering a function $f \in L_2(\Omega)$, we have the discretized truncated orthogonal projection $\mathcal{T}_n: L_2(\Omega) \to \mathbb{P}_{n}(\Omega)$
\begin{equation*}
    \mathcal{T}_n f:=\sum_{\ell=1}^{d_n} \hat{f}_{\ell} \Phi_{\ell}=\sum_{\ell=1}^{d_n} \langle f, \Phi_{\ell} \rangle \Phi_{\ell},
\end{equation*}
where $\{\hat{f}_{\ell} \}_{\ell=1}^{d_n}$ are the Fourier coefficients
\begin{equation}\label{equ:fourier}
    \hat{f}_{\ell} :=\langle f, \Phi_{\ell} \rangle = \int_{\Omega} f \Phi_{\ell} {\rm{d}}\omega, \qquad \forall \ell = 1, \cdots, d_n.
\end{equation}

To approximate the Fourier coefficients in \eqref{equ:fourier}, we use an exact quadrature formula for $\mathbb{P}_{2n}(\Omega)$ with nodes $\mathcal{X}_N=\{\mathbf{x}_1,\ldots, \mathbf{x}_{N}  \} \subset \Omega$ and positive weights $\mathbf{w}=\{w_1,\ldots, w_N\}$,
\begin{equation}\label{equ:quadrature}
    \int_{\Omega} p(\mathbf{x}) {\rm{d}}\omega = \sum_{j=1}^{N} w_j p(\mathbf{x}_j) \qquad \forall p \in \mathbb{P}_{2n}(\Omega).
\end{equation}
Corresponding to the inner product in \eqref{equ:inner}, Sloan introduced the ``discrete (semi) inner product''
\begin{equation}\label{equ:semiinner}
    \langle f, g\rangle_N = \sum_{j=1}^{N} w_jf(\mathbf{x}_j)g(\mathbf{x}_j) \qquad \forall f,g \in \mathcal{C}(\Omega),
\end{equation}
in which the exact integral is replaced by the quadrature rule. Based on the discrete (semi) inner product with quadrature exactness $2n$, Sloan originally proposed in \cite{sloan1995hyperinterpolation} the \emph{hyperinterpolation} $\mathcal{L}_n f$  which is a discretization of the $L_2$ orthogonal projection $\mathcal{T}_n f$ of $f \in \mathcal{C}(\Omega)$, i.e., $\mathcal{L}_n: \mathcal{C}(\Omega) \to \mathbb{P}_n(\Omega) $ as
\begin{equation}\label{equ:hyperinterpolation}
    \mathcal{L}_n f :=\sum_{\ell=1}^{d_n} \langle f, \Phi_{\ell} \rangle_{N} \Phi_{\ell}.
\end{equation}
For every $f \in \mathcal{C}(\Omega)$, hyperinterpolant $\mathcal{L}_{n}f$ satisfies the basic estimate
\begin{equation}\label{equ:error}
\| \mathcal{L}_{n}f - f\|_2 \leq 2V^{1/2}E_n(f) ,
\end{equation}
where $E_n(f):=\inf \{ \| f-p\|_{\infty};p \in \mathbb{P}_{n}(\Omega) \}$ and $E_n(f)$ tends to 0 as $n$ approaches infinity.

Sloan revealed that $\mathcal{L}_nf$ is the best discrete least squares approximation (weighted by quadrature weights) of $f$ at the quadrature points in \cite{sloan1995hyperinterpolation}. Consider the following discrete weighted least squares approximation problem
\begin{equation}\label{equ:approximationproblem}
\min\limits_{p\in\mathbb{P}_n(\Omega)}~~\left\{\sum_{j=1}^{N} w_j[p(\mathbf{x}_j)-f(\mathbf{x}_j)]^2 \right\}
\end{equation}
with $p(\mathbf{x})= \sum_{\ell=1}^{d_n} \alpha_{\ell}\Phi_{\ell}(\mathbf{x}) \in \mathbb{P}_{n}(\Omega)$, or equivalently
\begin{equation*}%\label{equ:optimizationproblem}
\min\limits_{\bm{\alpha}\in\mathbb{R}^{d_n}}~~ \|\mathbf{W}^{1/2}(\mathbf{A}\bm{\alpha}-\mathbf{f})\|_2^2,
\end{equation*}
where $\mathbf{W}=\text{diag}(w_1,\ldots,w_N) \in \mathbb{R}^{N \times N}$, $\mathbf{A} \in\mathbb{R}^{N\times d_n}$ with $[\mathbf{A}]_{j\ell}=\Phi_{\ell}(\mathbf{x}_j)$, $\bm{\alpha}=[\alpha_1,\ldots,\alpha_{d_n}]^{\rm{T}}\in\mathbb{R}^{d_n}$ and $\mathbf{f}=[f(\mathbf{x}_1),\ldots,f(\mathbf{x}_N)]^{\rm{T}}\in\mathbb{R}^N$ are two column vectors (recall $\mathbf{x}_j\in\mathbb{R}^s$).

Then we conclude the above discussion as the following result.
\begin{lemma}[{\cite[Lemma 5]{sloan1995hyperinterpolation}}]\label{prop:hyperinterpolation}
Given $f\in\mathcal{C}(\Omega)$, let $\mathcal{L}_n f\in\mathbb{P}_n(\Omega)$ be defined by \eqref{equ:hyperinterpolation}, where the quadrature exactness of the corresponding quadrature formula is $2n$. Then $\mathcal{L}_n f$ is the unique solution to the approximation problem \eqref{equ:approximationproblem}.
\end{lemma}

\begin{remark}
The quadrature exactness $2n$ in the approximation scheme of hyperinterpolation can be relaxed to $n+k$ with $0 < k \leq n$, which implies that the potential quadrature rules for constructing hyperinterpolation can be significantly enriched and the number of quadrature points can be considerably reduced  \cite{an2022exactness}. 
\end{remark}

\subsection{Lasso hyperinterpolation}\label{sec:lassohyper}
Loosely speaking, Lasso hyperinterpolation makes use of a \emph{soft thresholding operator} \cite{donoho1994ideal} 
\begin{equation*}
 \eta_{S}(a,k):=\max(0,a-k)+\min(0,a+k)   
\end{equation*}
to filter all hyperinterpolation coefficients $\langle{f,\Phi_{\ell}} \rangle_N$.

\begin{definition}[\cite{an2021lasso}]
Given a quadrature rule \eqref{equ:quadrature} with exactness $2n$, a \emph{Lasso hyperinterpolation} of $f$ onto $\mathbb{P}_{n}(\Omega)$ is defined as
\begin{equation}\label{equ:lassohyperinterpolation}
\mathcal{L}_n^{\lambda}{f}:=\sum_{\ell=1}^{d_n} \eta_{S}(\langle{f,\Phi_{\ell} \rangle_N } ,\lambda \mu_{\ell} )\Phi_{\ell},
\end{equation}
where $\lambda>0$ is the regularization parameter and $\{\mu_{\ell}\}_{\ell=1}^{d_n}$ is a set of positive penalty parameters.
\end{definition}
It has been proved in \cite{an2021lasso} that Lasso hyperinterpolation corresponds to an $\ell_1$-regularized least squares problem
\begin{equation}\label{equ:lasso}
\min\limits_{p_{\lambda}\in\mathbb{P}_n(\Omega)}~~\left\{\frac{1}{2} \sum_{j=1}^{N} w_j[p(\mathbf{x}_j)-f(\mathbf{x}_j)]^2 + \lambda \sum_{\ell=1}^{d_n} \mu_{\ell}|\gamma_{\ell}| \right\}
\end{equation}
with $p_{\lambda}(\mathbf{x})=\sum_{\ell=1}^{d_n} \gamma_{\ell}^{\lambda} \Phi_{\ell}(\mathbf{x}) \in \mathbb{P}_{n}(\Omega)$, or equivalently
\begin{equation*}%\label{equ:lassoMatrix}
    \min\limits_{\bm{\gamma}^{\lambda} \in \mathbb{R}^d} \:\:  \frac{1}{2} \| \mathbf{W}^{1/2}(\mathbf{A}\bm{\gamma}^{\lambda} - \mathbf{f})\|_2^2 + \lambda\|\mathbf{R}_1\bm{\gamma}^{\lambda}\|_1,
\end{equation*}
where $\lambda$ is a positive regularization parameter,$\bm{\gamma}^{\lambda}=[\gamma_1^{\lambda},\ldots,\gamma_{d_n}^{\lambda}]^{\text{T}}\in\mathbb{R}^{d_n}$ and $\mathbf{R}_1 = \text{diag}(\mu_1, \ldots, \mu_{d_n})$. In this paper, the subscript and superscript in $p_{\lambda}$ and $\bm{\gamma}^{\lambda}$ indicate that their specific forms are related to the regularization parameter $\lambda$. 

As shown in \cite{an2021lasso}, Lasso hyperinterpolation provides superior denoising performance compared to filtered hyperinterpolation \cite{sloan2012filtered}. Here, we aim to compare the denoising effectiveness of Lasso hyperinterpolation with that of hard thresholding hyperinterpolation.

\subsection{Hard thresholding hyperinterpolation}\label{sec:hardhyper}
Now we incorporate a hard thresholding operator into the hyperinterpolation operator $\mathcal{L}_n$. Based on the original idea of hard thresholding method proposed by Donoho and Johnstone \cite{donoho1994ideal}, where only a few wavelet coefficients contribute to the signal, we consider threshold rules that retain only observed data that exceed a multiple of the noise level. For the hyperinterpolant $\mathcal{L}_n f = \sum_{\ell=1}^{d_n}\langle f, \Phi_{\ell} \rangle_N \Phi_{\ell}$, this incorporation is achieved by solving the $\ell_0$-regularized weighted least squares problem
\begin{equation}\label{equ:hardthresholdapproximation}
\min\limits_{p_{\lambda}\in\mathbb{P}_n(\Omega)}~~\left\{\sum_{j=1}^{N} w_j[p_{\lambda}(\mathbf{x}_j)-f(\mathbf{x}_j)]^2
+ \lambda^2 \sum_{\ell=1}^{d_n} |\beta_{\ell}^{\lambda}|_0\right\}
\end{equation}
where $ p_{\lambda}(\mathbf{x})=\sum_{\ell=1}^{d_n}\beta_{\ell}^{\lambda} \Phi_{\ell}(\mathbf{x}) \in \mathbb{P}_{n}(\Omega)$, $\lambda>0$ is the \emph{regularization parameter} and $|\beta_{\ell}^{\lambda}|_{0}$ denotes the $\ell_0$-norm in one dimension, that is
\begin{equation*}
    \forall \beta_{\ell}^{\lambda} \in \mathbb{R},\quad |\beta_{\ell}^{\lambda}|_{0}:=
\left\{
\begin{array}{ll}
0, & \text{if}~~ \beta_{\ell}^{\lambda}=0, \\
1, & \text{if}~~ \beta_{\ell}^{\lambda}\ne0.
\end{array}
\right.
\end{equation*}
Problem \eqref{equ:hardthresholdapproximation} also amounts to the following
\begin{equation*}%\label{equ:hardthresholdoptimization}
    \min\limits_{\bm{\beta}^{\lambda} \in \mathbb{R}^d} \:\: \| \mathbf{W}^{1/2}(\mathbf{A}\bm{\beta}^{\lambda} - \mathbf{f})\|_2^2 + \lambda^2\|\bm{\beta}^{\lambda}\|_0, 
\end{equation*}
where $\lambda$ is a positive regularization parameter, $\bm{\beta}^{\lambda}=[\beta_1^{\lambda},\cdots,\beta_{d_n}^{\lambda}]^{\text{T}}\in\mathbb{R}^{d_n}$ and $\|\bm{\beta}^{\lambda}\|_0$ represents the number of nonzero elements of $\bm{\beta}^{\lambda}$. 

The reason $\lambda^2$ appears in equation \eqref{equ:hardthresholdapproximation} while $\lambda$ is used in equation \eqref{equ:lasso} is to set a consistent thresholding value for both Lasso and hard thresholding hyperinterpolations when selecting coefficients. Consequently, it is convenient for us to make numerical comparisons, see section \ref{sec:examples}.

Then we give definitions of \emph{hard thresholding operator} and \emph{hard thresholding hyperinterpolation}, respectively.

\begin{definition}[\cite{donoho1994ideal}]\label{def:hardthresholdingoperator}
The \emph{hard thresholding operator}, denoted by $\eta_{H}(a,k)$, is defined as
\[\eta_{H}(a,k):=\left\{\begin{array}{ll}
a,  & \text{if}~~|a|> k ,\\
0,  & \text{if}~~|a|\leq k.
\end{array}\right. \]
\end{definition}

\begin{definition}\label{def:hardthresholdhyperinterpolation}
Given a quadrature rule \eqref{equ:quadrature} with exactness $2n$, {a \emph{hard thresholding hyperinterpolation}} of $f$ onto $\mathbb{P}_n(\Omega)$ is defined as
\begin{equation}\label{equ:hardthresholdhyperinterpolation}
\mathcal{H}_{n}^{\lambda}{f}:=\sum_{\ell=1}^{d_n}  \eta_{H}(\langle{f,\Phi_{\ell}}\rangle_N ,\lambda)\Phi_{\ell},\qquad \lambda>0.
\end{equation}
\end{definition}

\begin{remark}
When $\lambda $ takes zero, hard thresholding hyperinterpolation operator $\mathcal{H}_{n}^{\lambda}$ becomes classical hyperinterpolation operator $\mathcal{L}_n$. 
\end{remark}

%\begin{remark}
%Hard threholding hyperinterpolation operator $\mathcal{H}_n^{\lambda}$ is nonlinear.
%\end{remark}

Now we show that the hard thresholding hyperinterpolation  $\mathcal{H}_n^{\lambda} f$ is indeed the unique solution to the $\ell_0$-regularized least squares problem \eqref{equ:hardthresholdapproximation}.

\begin{theorem}\label{thm:ExactSolution}
Let $\hardhyper{f} \in \mathbb{P}_{n}(\Omega)$ be defined by \eqref{equ:hardthresholdhyperinterpolation}, and adopt conditions of Lemma \ref{prop:hyperinterpolation}. Then $\hardhyper{f} \in \mathbb{P}_{n}(\Omega)$ is the unique solution to the regularized least squares approximation problem \eqref{equ:hardthresholdapproximation}.
\end{theorem}

\begin{proof}
Let $\bm{\alpha}=[\alpha_1, \ldots, \alpha_{d_n}]^{\rm{T}} \in \mathbb{R}^{d_n}$ with $\alpha_{\ell}=\langle f, \Phi_{\ell} \rangle_N$ for $\ell=1,\ldots, d_n$. 
Then the first-order condition is satisfied by taking the first derivative of $\normtwo{\mathbf{W}^{1/2}(\mathbf{A}\bm{\alpha}-\mathbf{f})}^2$ with respect to $\bm{\alpha}$ and setting it equal to zero
\begin{equation}\label{equ:firstorder}
2\mathbf{A}^{\rm{T}}\mathbf{WA}\bm{\alpha}-2\mathbf{A}^{\rm{T}}\mathbf{Wf}=\mathbf{0}.
\end{equation}
We assert that $\mathbf{A}^{\rm{T}}\mathbf{WA}$ is an identity matrix:
\begin{equation}\label{equ:identity}
[\mathbf{A}^{\rm{T}}\mathbf{WA}]_{ik} =\sum\limits_{j=1}^{N} w_j\Phi_{i}(\mathbf{x}_j)\Phi_{k}(\mathbf{x}_j)=\discreteinner{\Phi_i,\Phi_k}=\delta_{ik}, \quad 1\leq i,k \leq d_n.
\end{equation}
By \eqref{equ:firstorder} and \eqref{equ:identity}, we obtain $\bm{\alpha}=\mathbf{A}^{\rm{T}}\mathbf{Wf}$. Let $$F(\bm{\beta}^{\lambda}):=\normtwo{\mathbf{W}^{1/2}(\mathbf{A} \bm{\beta}^{\lambda} - \textbf{f})}^2 + \lambda^2 \normzero{ \bm{\beta}^{\lambda}}.$$ Then we can decompose $F(\bm{\beta}^{\lambda})$ as the following:
\begin{equation*}
   \begin{aligned}
F(\bm{\beta}^{\lambda}) 
&= (\bm{\beta}^{\lambda})^{\rm{T}} \mathbf{A}^{\rm{T}}\mathbf{WA} \bm{\beta}^{\lambda} - 2 (\bm{\beta}^{\lambda})^{\rm{T}}\mathbf{A}^{\rm{T}}\mathbf{Wf} + \textbf{f}^{\rm{T}} \mathbf{Wf} + \lambda^2 \normzero{ \bm{\beta}^{\lambda}} ,\\
&=(\bm{\beta}^{\lambda})^{T}\bm{\beta}^{\lambda} - 2(\bm{\beta}^{\lambda})^{\rm{T}} \bm{\alpha} +\bm{\alpha}^{\rm{T}}\bm{\alpha} - \bm{\alpha}^{\rm{T}}\bm{\alpha} + \textbf{f}^{\rm{T}} \mathbf{Wf}  +  \lambda^2 \normzero{ \bm{\beta}^{\lambda}},\\
&= \textbf{f}^{\rm{T}} \mathbf{Wf}-\bm{\alpha}^{\rm{T}}\bm{\alpha} + \normtwo{\bm{\beta}^{\lambda} - \bm{\alpha}}^2 + \lambda^2 \| \bm{\beta}\|_0,\\
&=C+ \sum\limits_{\ell =1}^{d_n}\squarebrackets{  (\beta_{\ell}^{\lambda} - \alpha_{\ell})^{2}  + \lambda^2 |\beta_{\ell}^{\lambda}|_0},
\end{aligned} 
\end{equation*}
where $C=\sum_{j=1}^{N} w_jf^2(\textbf{x}_j)- \sum_{\ell=1}^{d_n}\alpha_{\ell}^2$. Thus, solving problem \eqref{equ:hardthresholdapproximation} is equivalent to solve $d_n$ independent one-dimensional problems, i.e.,
\begin{equation*}
 g(\beta_{\ell}^{\lambda}):=\left\{\begin{array}{ll}
\alpha_{\ell}^2 +r(\beta_{\ell}^{\lambda}) , &\beta_{\ell}^{\lambda} \neq 0, \\
\alpha_{\ell}^2, & \beta_{\ell}^{\lambda}=0,
\end{array}
\right.   
\end{equation*}
where $r(\beta_{\ell}^{\lambda}):=(\beta_{\ell}^{\lambda})^2 - 2\beta_{\ell}^{\lambda}\alpha_{\ell}+ \lambda^2$. Since $r(\beta_{\ell}^{\lambda})$ is quadratic with respect to the variable $\beta_{\ell}^{\lambda}$, we find $\min r(\beta_{\ell}^{\lambda}) = r(\alpha_{\ell})$ and have the discriminant 
\begin{equation*}
    \Delta :=(2\alpha_{\ell})^2-4\cdot\lambda^2 = 4\alpha_{\ell}^2 -4\lambda^2.
\end{equation*}
There are two cases that we need to consider:
\begin{enumerate}
\item[(i)] If $r(\alpha_{\ell}) <0$,
then $\Delta > 0$ which means that $|\alpha_{\ell}|>\lambda$. We deduce that $g(\beta_{\ell}^{\lambda})$ arrives at its minimum when $\beta_{\ell}^{\lambda} = \alpha_{\ell}  ~~(|\alpha_{\ell}| > \lambda)$;
\item[(ii)] If $r(\alpha_{\ell}) \geq 0$, then $\Delta \leq 0$ which means that $|\alpha_{\ell}|\leq \lambda$. In this case, we obtain that
\begin{equation*}
    \min g(\beta_{\ell}^{\lambda}) = g(0) \leq  g(\alpha_{\ell}) ,
\end{equation*}
where the inequality becomes equality when $|\alpha_{\ell}|=\lambda$.
\end{enumerate}
Therefore, for $\ell=1,\ldots,d_n$, we obtain
\begin{equation*}
 \beta_{\ell}^{\lambda}=\hardoperatornew{\alpha_{\ell}} = \left\{\begin{array}{ll}
\alpha_{\ell},  & \text{if}~~|\alpha_{\ell}|> \lambda ,\\
0,  & \text{if}~~|\alpha_{\ell}|\leq \lambda.
\end{array}\right.    
\end{equation*}

Recall that $\alpha_{\ell}=\discreteinner{f, \Phi_{\ell}}$ for all $\ell=1,\ldots, d_n$. We deduce that the polynomial constructed with coefficients $\beta_{\ell}^{\lambda}=\hardoperatornew{\alpha_{\ell}},\ell=1,\ldots, d_n$, is indeed $\hardhyper{f} $.  
\end{proof}

\section{Characterizations of hard thresholding hyperinterpolation}\label{sec:HTH}
In this section, we explore some algebraic and geometric properties of hard thresholding hyperinterpolation, respectively.

\subsection{Algebraic properties}
First, we investigate the idempotent property of hard thresholding hyperinterpolation operator $\mathcal{H}_n^{\lambda}$, and the commutative law between $\mathcal{H}_n^{\lambda}$ and hyperinterpolation operator $\mathcal{L}_n$.
\begin{theorem}\label{thm:orthogonalprojection}
Let $ \hyper{f}, \hardhyper{f}  \in \mathbb{P}_{n}(\Omega)$ be defined by \eqref{equ:hyperinterpolation} and \eqref{equ:hardthresholdhyperinterpolation}, respectively, and adopt conditions of Lemma \ref{prop:hyperinterpolation}. Then
\begin{itemize}
 \item[(a)] $\mathcal{H}_n^{\lambda}$ is idempotent: $\hardhyper(\hardhyper{f})=\hardhyper{f}$,
 \item[(b)] $\hardhyper(\hyper{f})=\hyper(\hardhyper{f})=\hardhyper{f}$.
\end{itemize}
\end{theorem}

\begin{proof}
(a)  Since $\{\Phi_{1},\ldots,\Phi_{d_n}\}$ is an orthonormal basis and by \eqref{equ:hardthresholdhyperinterpolation}, we have
\begin{equation*}
\begin{aligned}
\hardhyper(\hardhyper{f}) &= \sum\limits_{\ell=1}^{d_n} \hardoperator{ \sum\limits_{j=1}^{d_n} \hardoperator{f,\Phi_j}\Phi_{j}  , \Phi_{\ell}  }\Phi_{\ell}, \\
&=\sum\limits_{\ell=1}^{d_n} \hardoperator{\hardoperator{f,\Phi _{\ell}}\Phi_{\ell}, \Phi_{\ell}}\Phi_{\ell} ,\\
&=\sum\limits_{\ell=1}^{d_n} \hardoperatornew{\hardoperator{f,\Phi_{\ell}}}\Phi_{\ell} .\\
\end{aligned}
\end{equation*}
There are two cases that we need to consider:
\begin{enumerate}
\item[(i)] If $|\discreteinner{f,\Phi_{\ell}}| > \lambda $, then \[\hardoperatornew{\hardoperator{f, \Phi_{\ell}}} =\hardoperator{f, \Phi_{\ell}}= \discreteinner{f,\Phi_{\ell}}.\]
\item[(ii)] If $|\discreteinner{f,\Phi_{\ell}}|\leq \lambda$, then
\[\hardoperatornew{\hardoperator{f, \Phi_{\ell}}} =\hardoperatornew{0}= 0.\]
\end{enumerate}
Combing the above, we obtain $\hardhyper(\hardhyper{f})=\hardhyper{f}$.

(b) Recall that $\hyper{f},\hardhyper{f} \in \mathbb{P}_n(\Omega)$ defined by \eqref{equ:hyperinterpolation} and \eqref{equ:hardthresholdhyperinterpolation}, respectively, we obtain
\begin{equation*}
\begin{aligned}
\hyper(\hardhyper{f}) & = \sum\limits_{j=1}^{d_n} \discreteinner{ \sum\limits_{\ell=1}^{d_n} \hardoperator{f,\Phi_{\ell}}\Phi_{\ell}, \Phi_{j} }\Phi_{j} , \\
& = \sum\limits_{j=1}^{d_n} \discreteinner{  \hardoperator{f,\Phi_{j}}\Phi_{j}, \Phi_{j} }\Phi_{j} , \\
& =  \sum\limits_{j=1}^{d_n}  \hardoperator{f,\Phi_{j}} \Phi_{j},
\end{aligned}
\end{equation*}
and
\begin{equation*}
\begin{aligned}
\hardhyper(\hyper{f}) & = \sum\limits_{\ell=1}^{d_n} \hardoperator{ \sum\limits_{j=1}^{d_n}\discreteinner{f,\Phi_j}\Phi_j, \Phi_{\ell}}\Phi_{\ell} , \\
&= \sum\limits_{\ell=1}^{d_n} \hardoperator{  \discreteinner{f,\Phi_{\ell}}\Phi_{\ell}, \Phi_{\ell}}\Phi_{\ell} , \\
&= \sum\limits_{\ell=1}^{d_n} \hardoperator{ f,\Phi_{\ell} }\Phi_{\ell} .
\end{aligned}
\end{equation*}
Thus, we have completed the proof.  
\end{proof}

\begin{remark}
Notice that hard thresholding hyperinterpolation operator is not symmetric in the sense of
\begin{equation*}
    \langle  \mathcal{H}_{n}^{\lambda} f , g \rangle_N \neq \langle  f ,\mathcal{H}_{n}^{\lambda} g \rangle_N, 
\end{equation*}
where $f,g \in \mathcal{C}(\Omega)$. 
\end{remark}

\subsection{Geometric properties}
The following lemma describes the geometric properties of $\hardhyper$, which is very helpful in proving Theorem \ref{thm:Noise}. 
\begin{lemma}\label{lem:mainlemma}
Under conditions of Theorem \ref{thm:ExactSolution},
\begin{itemize}
  \item[(a)]$\discreteinner{f-\hardhyper{f},\hardhyper{f}}=0$,
  \item[(b)] $\discreteinner{\hardhyper{f},\hardhyper{f}}+\discreteinner{f-\hardhyper{f},f-\hardhyper{f}}=
\discreteinner{f,f}$,
  \item[(c)]$\discreteinner{\hardhyper{f},\hardhyper{f}}\leq\discreteinner{f,f}$.
\end{itemize}
\end{lemma}

\begin{proof}
(a) Let $\alpha_{\ell}=\discreteinner{f,\Phi_{\ell}},\: \ell=1,\ldots,d_n$. Then we obtain
\begin{equation*}
\begin{aligned}
\discreteinner{f-\hardhyper{f},\hardhyper{f}}&=\discreteinner{f- \sum\limits_{\ell=1}^{d_n} \hardoperatornew{\alpha_{\ell}}\Phi_{\ell}, \sum\limits_{k=1}^{d_n} \hardoperatornew{\alpha_{k}}\Phi_{k}}, \\
&=\sum\limits_{k=1}^{d_n}\hardoperatornew{\alpha_k}
\discreteinner{f-\sum\limits_{\ell=1}^{d_n}\hardoperatornew{\alpha_{\ell}}\Phi_{\ell},\Phi_k},
\end{aligned}
\end{equation*}
and
\begin{equation*}
\begin{aligned}
\discreteinner{f-\sum\limits_{\ell=1}^{d_n}\hardoperatornew{\alpha_{\ell}}\Phi_{\ell},\Phi_k}&=\discreteinner{f,\Phi_k}-\discreteinner{\sum_{\ell=1}^{d_n}\hardoperatornew{\alpha_{\ell}}\Phi_{\ell},\Phi_k},\\
&=\alpha_{k}-\hardoperatornew{\alpha_{k}}.
\end{aligned}
\end{equation*}
The fact
\begin{equation*}
\hardoperatornew{\alpha_{k}} = \left\{\begin{array}{ll}
\alpha_{k},  & \text{if}~~|\alpha_{k}|> \lambda ,\\
0,  & \text{if}~~|\alpha_{k}|\leq \lambda,
\end{array}\right.    
\end{equation*}
gives the following
\begin{equation*}
     \hardoperatornew{\alpha_k} (\alpha_{k}-\hardoperatornew{\alpha_{k}})  = ( \hardoperatornew{\alpha_k}\alpha_{k}-(\hardoperatornew{\alpha_k})^{2})=0.
\end{equation*}
Therefore, we obtain $\discreteinner{f-\hardhyper{f},\hardhyper{f}}=0$.

(b) By $\discreteinner{\hardhyper{f},\hardhyper{f}}=\discreteinner{f,\hardhyper{f}}$ and
\begin{equation*}
    \discreteinner{f-\hardhyper{f},f-\hardhyper{f}} = \discreteinner{f,f} -2\discreteinner{f,\hardhyper{f}} +\discreteinner{\hardhyper{f},\hardhyper{f}},
\end{equation*}
we obtain
\begin{equation}\label{equ:mainlemma2}
 \discreteinner{\hardhyper{f},\hardhyper{f}}  + \discreteinner{f-\hardhyper{f},f-\hardhyper{f}} = \discreteinner{f,f}.
\end{equation}

(c) From  \eqref{equ:mainlemma2} and the positiveness of $\discreteinner{f-\hardhyper{f},f-\hardhyper{f}}$, we immediately have
\begin{equation*}
    \discreteinner{\hardhyper{f},\hardhyper{f}} \leq \discreteinner{f,f} .
\end{equation*}
Thus, we have completed the proof. 
\end{proof}

\begin{figure}[htbp]
  \centering
  \includegraphics[scale=0.9,clip]{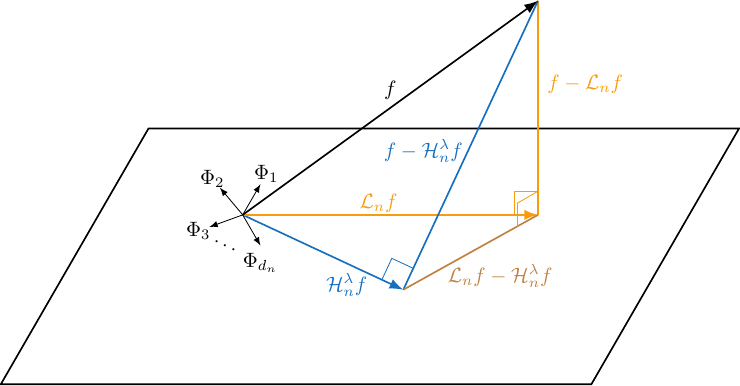}\\
  \caption{The geometric interpretation of hyperinterpolation and hard thresholding hyperinterpolation satisfying the Pythagorean theorem with respect to the discrete (semi) inner product \eqref{equ:semiinner}.}\label{Pythagorean}
\end{figure}

\begin{remark}\label{rema1}
Recall that for hyperinterpolation, we have
\begin{equation*}
\langle f-\mathcal{L}_n f, p \rangle_N=0 
\end{equation*}
for any $p \in\mathbb{P}_{n}(\Omega)$, and
\begin{equation*}
 \langle \mathcal{L}_n f, \mathcal{L}_nf \rangle_N +\langle f-\mathcal{L}_n f, f- \mathcal{L}_n f \rangle_N = \langle f, f \rangle_N.
\end{equation*}
Both classical hyperinterpolation and hard thresholding hyperinterpolation satisfy the Pythagorean theorem with respect to the discrete (semi) inner product \eqref{equ:semiinner}, as shown in Figure \ref{Pythagorean}. However, Lasso hyperinterpolation does not possess this geometric property \cite{an2021lasso}. 
\end{remark}

\begin{proposition}\label{prop:BasisNotInvariant}
Let $\{q_i\}_{i=1}^{d_n}$ and $\{\Phi_i\}_{i=1}^{d_n}$ be two orthonormal bases of $\mathbb{P}_n(\Omega)$ and $\Phi_{i}=\sum_{j=1}^{d_n}a_{ij}q_{j}$ with $a_{ij}=\langle q_i, \Phi_{j} \rangle_N$, and let $ b_{ij}= \sum_{k=1}^{d_n}a_{ik}a_{kj}$ for $i,j=1,\ldots,d_n$. Then hard thresholding hyperinterpolation is not invariant under a change of basis. That is, 
\begin{equation*}%\label{equ:BasisNotInvariant}
    \hardhyper{f} \neq \sum_{i,j=1}^{d_n}\hardoperator{f,q_i}b_{ij}q_j.
\end{equation*}
\end{proposition}

\begin{proof}
We assert that hard thresholding hyperinterpolation is not basis invariant, i.e.,
\begin{equation*}%\label{equ:Invariant1}
    \begin{aligned}
        \mathcal{H}_n^{\lambda} f &=\sum_{i=1}^{d_n} \eta_{H}(\langle f, \Phi_i \rangle_N, \lambda) \Phi_i , \\&= \sum_{i=1}^{d_n} \eta_{H}\left(\left\langle f, \sum_{j=1}^{d_n} a_{ij} q_j \right\rangle_N, \:\lambda \right) \sum_{j=1}^{d_n} a_{ij} q_j ,\\
        &\neq \sum_{i=1}^{d_n} \eta_{H}(\langle f, q_i \rangle_N, \lambda) \sum_{j=1}^{d_n}a_{ij} \cdot\sum_{j=1}^{d_n} a_{ij}q_j ,\\
        & = \sum_{i,j=1}^{d_n} \eta_{H}(\langle f, q_i \rangle_N,\lambda) b_{ij}q_j,
    \end{aligned}
\end{equation*}
where the inequality holds because of the presence of a hard threhoslding operator.   
\end{proof}

\section{Error analysis}\label{sec:erroranalysis}
\subsection{Error analysis in the \texorpdfstring{$L_2(\Omega)$}{}-norm}\label{sec:L2NormError}
In the following, the norms which appear in the theorems are defined by
\begin{equation*}
\begin{aligned}
\|g\|_2 &:= \left( \int_\Omega g^2 \text{d}\omega\right)^{1/2}, \quad g \in L_2(\Omega),\\  
\|g\|_{\infty} &:= \sup\limits_{\mathbf{x} \in \Omega} |g(\mathbf{x})|,  \qquad \: \:\:g \in \mathcal{C}(\Omega).  
\end{aligned}
\end{equation*}
For any operator $\mathcal{U}_n: \mathcal{C}(\Omega) \to L_2(\Omega)$, its operator norm is defined as
\begin{equation*}
    \|\mathcal{U}_n\|_{\text{op}} :=\sup_{g\in \mathcal{C}(\Omega), g \neq 0} \frac{\|\mathcal{U}_n g\|_{2}}{\|g\|_{\infty}}.
\end{equation*}

The error of the best approximation to $f\in \mathcal{C}(\Omega)$ by polynomials of degree at most $n$ is defined by 
\begin{equation*}
    E_n(f)=\inf\{\|f-p\|_{\infty};p \in \mathbb{P}_n(\Omega)\}.
\end{equation*}

Note that there are various techniques available for stable approximation, including piecewise interpolation, we have chosen to maintain our focus on hyperinterpolation and its variants due to their unique advantages in providing stable, high-order polynomial approximations over compact manifolds. This is evidenced by the operator norm estimations in  \cite{MR2274179} and  \cite{sloan2012filtered}, and notably, as shown in \cite{sloan1997interpolation}, the $L_2$ operator norm of hyperinterpolation is strictly less than that of interpolation over the unit sphere.

The following theorem concentrates on the operator norms of Lasso hyperinterpolation operator $\mathcal{L}_n^{\lambda}$, hard thresholding  hyperinterpolation operator $\mathcal{H}_n^{\lambda}$, and hyperinterpolation operator $\mathcal{L}_n$. 

\begin{theorem}\label{thm:lassohardclassicalhyper}
Adopt conditions of Lemma \ref{prop:hyperinterpolation}.
In Lasso hyperinterpolation \eqref{equ:lassohyperinterpolation} , set all $\{\mu_{\ell}\}_{\ell=1}^{d_n}$ to 1.  
For any given $\lambda \in[0,\,\norminfty{\mathbf{A}^{{\rm{T}}}\mathbf{Wf}}) $,  we have
\begin{equation}
\normtwo{\lassohyper{f}} \leq \normtwo{\hardhyper{f}} \leq \normtwo{\hyper{f}},
\end{equation}
where the equality holds when $\lambda=0$. Moreover,
\begin{equation}\label{equ:operator}
    \| \mathcal{L}_n^{\lambda} \|_{\rm{op}} \leq \| \mathcal{H}_n^{\lambda} \|_{\rm{op}} \leq \| \mathcal{L}_n \|_{\rm{op}}.
\end{equation}
\end{theorem}

\begin{proof}
Let $\alpha_{\ell}:=\discreteinner{f,\Phi_{\ell}}$ for all $1 \leq \ell \leq d_n$. Since $\{\Phi_{1},\ldots,\Phi_{d_n}\}$ is an orthonormal basis of $\mathbb{P}_{n}(\Omega)$ and by Parseval's identity, we have
\begin{equation*}
    \normtwo{\hyper{f}} = \brackets{\sum\limits_{\ell=1}^{d_n} \absolute{\alpha_{\ell}}^2}^{\frac{1}{2}} \quad \text{and} \quad \normtwo{\hardhyper{f}} = \brackets{\sum\limits_{\ell=1, \alpha_{\ell} > \lambda }^{d_n} \alpha_{\ell}^2 + \sum\limits_{\ell=1, \alpha_{\ell} < -\lambda }^{d_n} \alpha_{\ell}^2}^{\frac{1}{2}},
\end{equation*}
and under the condition $\{\mu_{\ell}\}_{\ell=1}^{d_n}$ all being 1, 
\begin{equation*}
\begin{aligned}
\normtwo{\lassohyper{f}} &= \brackets{\sum\limits_{\ell=1, \alpha_{\ell} > \lambda\mu_{\ell} }^{d_n} (\alpha_{\ell}-\lambda\mu_{\ell})^2 + \sum\limits_{\ell=1, \alpha_{\ell} < -\lambda\mu_{\ell} }^{d_n} (\alpha_{\ell}+\lambda\mu_{\ell})^2 }^{\frac{1}{2}}, \\
&= \brackets{\sum\limits_{\ell=1, \alpha_{\ell} > \lambda }^{d_n} (\alpha_{\ell}-\lambda)^2 + \sum\limits_{\ell=1, \alpha_{\ell} < -\lambda }^{d_n} (\alpha_{\ell}+\lambda)^2 }^{\frac{1}{2}}.
\end{aligned}
\end{equation*}
If $\alpha_{\ell} > \lambda$, then $\alpha_{\ell}^{2} > (\alpha_{\ell}-\lambda)^2$; if $\alpha_{\ell} < -\lambda$, then $\alpha_{\ell}^2 > (\alpha_{\ell}+\lambda)^2$. Thus, it is clear that
\begin{equation*}
    \normtwo{\lassohyper{f}} \leq \normtwo{\hardhyper{f}} \leq \normtwo{\hyper{f}} .
\end{equation*}
When $\lambda=0$, the conclusion is obvious since both $\lassohyper$ and $\hardhyper$ become as $\hyper$. By the definition of the operator norm, we finish the proof.  
\end{proof}

Next, we examine the $L_2$ errors for hard threhoslding hyperinterpolation in approximating the test function given a specific noise $\epsilon$. Let $\{f^{\epsilon}(\mathbf{x}_j)\}_{j=1}^N$ represent the noisy observations of $\{f(\mathbf{x}_j)\}_{j=1}^N$ , where $f^{\epsilon}(\mathbf{x}_j) = f(\mathbf{x}_j) + \epsilon_j$ for $j=1,\ldots, N$. 

\begin{theorem}\label{thm:Noise}
Adopt conditions of Lemma \ref{thm:ExactSolution} and assume  $f^{\epsilon} = f+ \epsilon$ is a noisy version of $f$, where $\epsilon$ is some noise.
Then
\begin{equation}\label{equ:OperatorNorm_HardHyper}
    \|\mathcal{H}_n^{\lambda}f\|_2 \leq V^{1/2}\|f\|_{\infty},
\end{equation}
and
\begin{equation}\label{equ:L2errorNoise}
\|\mathcal{H}_n^{\lambda}f^{\epsilon} -f \|_2 \leq \| \mathcal{H}_n^{\lambda}f^{\epsilon} - \mathcal{L}_nf ^{\epsilon}\|_{2}  + 2V^{1/2}E_n(f) + V^{1/2} \|{\epsilon} \|_{\infty}.
\end{equation}
Thus
\begin{equation*}
\lim_{n \to \infty} \lim_{\lambda \to 0} \|\mathcal{H}_n^{\lambda}f ^{\epsilon}- f\|_2 =V^{1/2} \|{\epsilon} \|_{\infty} .
\end{equation*}
\end{theorem}

\begin{proof}
For any $f \in \mathcal{C}(\Omega)$, it is clear that $\mathcal{H}_n^{\lambda} f \in \mathbb{P}_{n}(\Omega)$. Then we have
\begin{equation*}
    \begin{aligned}
        \|\mathcal{H}_n^{\lambda}f\|_2^2 &=\langle \mathcal{H}_n^{\lambda}f , \mathcal{H}_n^{\lambda}f \rangle=\langle \mathcal{H}_n^{\lambda}f , \mathcal{H}_n^{\lambda}f \rangle_N \leq \langle f , f \rangle_N , \\
        & = \sum_{j=1}^{N} w_j (f(\mathbf{x}_j))^2 \leq \|f \|^{2}_{\infty} \sum_{j=1}^{N}w_j = V \|f \|^{2}_{\infty},
    \end{aligned}
\end{equation*}
where the first inequality follows from Lemma \ref{lem:mainlemma} (c). Consequently, we deduce:
\begin{equation*}
    \|\mathcal{H}_n^{\lambda}f\|_2 \leq V^{1/2} \|f \|_{\infty}.
\end{equation*}

Recalling that $\mathcal{L}_n p = p$ for any $p \in \mathbb{P}_n(\Omega)$ and letting $p^{\ast}$ be the best approximation of $f$ in $\mathbb{P}_n(\Omega)$, we have
\begin{equation*}
\begin{aligned}
\| \mathcal{L}_n f^{\epsilon} -f \|_2 &=\| \mathcal{L}_n (f^{\epsilon} - p^{\ast}) + p^{\ast} - f \|_2, \\
&\leq V^{1/2}\|f^{\epsilon} - p^{\ast} \|_{\infty} + V^{1/2}\|f - p^{\ast} \|_{\infty}, \\
&\leq V^{1/2}\|f^{\epsilon} - f \|_{\infty} +2 V^{1/2}E_{n}(f),
\end{aligned}
\end{equation*}
where the the second term on the right side in the second row follows from the Cauchy-Schwarz inequality. 

By direct computation, we obtain
\begin{equation*}
\begin{aligned}
\|\mathcal{H}_n^{\lambda}f^{\epsilon} -f \|_2 &\leq \| \mathcal{H}_n^{\lambda}f^{\epsilon} - \mathcal{L}_nf ^{\epsilon}\|_{2} + \| \mathcal{L}_n f^{\epsilon} -f \|_2, \\
& \leq \| \mathcal{H}_n^{\lambda}f^{\epsilon} - \mathcal{L}_nf ^{\epsilon}\|_{2}  + 2 V^{1/2}E_{n}(f) +  V^{1/2}\|{\epsilon} \|_{\infty}.
\end{aligned}
\end{equation*}
%Noting that $\mathcal{H}_n^{\lambda}f^{\epsilon} - \mathcal{L}_nf^{\epsilon} \in \mathbb{P}_n(\Omega)$ and using Parseval's identity,  we have 
%\begin{equation*}
%    \|\mathcal{H}_n^{\lambda}f^{\epsilon} - \mathcal{L}_nf ^{\epsilon}\|_2 =\left(\sum_{\ell=1,|\langle f^{\epsilon}, \Phi_{\ell}\rangle_N|\leq \lambda}^{d_n} |\langle f^{\epsilon}, \Phi_{\ell}\rangle_N|^2\right)^{\frac{1}{2}}=A_{n,\lambda,f^{\epsilon}}.
%\end{equation*}
Since 
\begin{equation*}
\lim_{\lambda \to 0} \|\mathcal{H}_n^{\lambda}f^{\epsilon} - \mathcal{L}_n f^{\epsilon}\|_2 =0   \quad \text{and} \quad  \lim_{n \to \infty} E_{n}(f) =0,
\end{equation*}
we have
\begin{equation*}
\lim_{n \to \infty} \lim_{\lambda \to 0} \|\mathcal{H}_n^{\lambda}f ^{\epsilon}- f\|_2 =V^{1/2} \|{\epsilon} \|_{\infty}.
\end{equation*}
Thus, we have completed the proof. 
\end{proof}

%\begin{remark}
%For Lasso hyperinterpolation, we obtain
%\begin{equation}
%\|\mathcal{L}_n^{\lambda}f^{\epsilon} -f \|_2 \leq \| \mathcal{L}_n^{\lambda}f^{\epsilon} - \mathcal{L}_nf ^{\epsilon}\|_{2}  + 2V^{1/2}E_n(f) + V^{1/2} \|f^{\epsilon} - f \|_{\infty}.
%\end{equation}
%Notably, we have
%\begin{equation*}
%\left\{
%    \begin{array}{l}
%\| \mathcal{H}_n^{\lambda}f^{\epsilon} - \mathcal{L}_nf ^{\epsilon}\|_{2}^2 = \sum\limits_{\ell=1, |\langle f^{\epsilon}, \Phi_{\ell} \rangle_N| \leq \lambda}^{d_n} \langle f^{\epsilon}, \Phi_{\ell} \rangle_N^2, \\
%\\
%\| \mathcal{L}_n^{\lambda}f^{\epsilon} - \mathcal{L}_nf ^{\epsilon}\|_{2}^2 = \sum\limits_{\ell=1, |\langle f^{\epsilon}, \Phi_{\ell} \rangle_N| \leq \lambda\mu_{\ell}}^{d_n} \langle f^{\epsilon}, \Phi_{\ell} \rangle_N^2 + \sum\limits_{\ell=1, |\langle f^{\epsilon}, \Phi_{\ell} \rangle_N| > \lambda\mu_{\ell}}^{d_n} (\lambda \mu_{\ell})^2.
%    \end{array}
%\right.
%\end{equation*}
%If we set $\{\mu_{\ell}\}_{\ell=1}^{d_n}$ all being 1, then the upper bound of the $L_2$ errors for hard thresholding hyperinterpolation is lower than that for Lasso hyperinterpolation. Compared to hard thresholding hyperinterpolation, Lasso introduces an additional term in the $L_2$ error related to the regularization parameter $\lambda$ and penalty parameters $\{ \mu_{\ell}\}_{\ell=1}^{d_n}$.  We will give a further quantitative analysis on the expansion of $L_2$ errors in detail in section \ref{sec:computational}.
%\end{remark}

\subsection{Error analysis in the uniform norm}
Recall that the function
\begin{equation}\label{equ:Christoffel}
    G_{n}(\mathbf{x},\mathbf{x}) = \sum_{\ell=1}^{d_n} \Phi_{\ell}^2(\mathbf{x})
\end{equation}
is known as the reciprocal of \emph{Christoffel function} {\cite{Nevai1986Christoffel}} associated with the measure $\text{d}\omega$ on $\Omega$. The reciprocal of \emph{Christoffel function} plays an important role in proving the uniform operator norms for $\|\mathcal{L}_n\|_{\infty}$ and $\|\mathcal{H}_n^{\lambda}\|_{\infty}$. Motivated by \cite[Proposition 1.1]{sommariva2014multivariate} and arguments in \cite{sommariva2021spheritri}, the error estimate \eqref{equ:error} for hyperinterpolation in $L_2$ norms can be extended to a corresponding version in the uniform norm, as detailed below.
\begin{lemma}[\cite{sommariva2014multivariate,sommariva2021spheritri}]\label{lem:hyperL2UNI}
Adopt conditions of Lemma \ref{prop:hyperinterpolation}. Let $C_n= \sqrt{\max\limits_{\mathbf{x} \in \Omega} G_{n}(\mathbf{x},\mathbf{x})}$ with $G_n$ defined by \eqref{equ:Christoffel}. Then
\begin{equation}\label{equ:inf_hyper}
    \| \mathcal{L}_n f\|_{\infty}  \leq C_n V^{1/2} \| f\|_{\infty}
\end{equation}
and
\begin{equation}\label{equ:inf_hyper_err}
    \| \mathcal{L}_n f - f\|_{\infty} \leq  (1+C_{n}V^{1/2})E_{n} (f).
\end{equation}
\end{lemma}

\begin{proof}
For any nonzero polynomial $p=\sum_{\ell=1}^{d_n} \langle p, \Phi_{\ell} \rangle_N \Phi_{\ell} \in \mathbb{P}_n({\Omega})$, we have
\begin{equation*}
\begin{aligned}
|p(\mathbf{x})|&=\left|\sum_{\ell=1}^{d_n} \langle p, \Phi_{\ell} \rangle_N \Phi_{\ell} (\mathbf{x})\right| \leq   \left(\sum_{\ell=1}^{d_n} \langle p, \Phi_{\ell} \rangle^2_N \right)^{1/2} (G_n(\mathbf{x},\mathbf{x}))^{1/2} , \\
 &\leq  \|p\|_2\left(\max_{\mathbf{x} \in \Omega} G_n(\mathbf{x},\mathbf{x})\right)^{1/2},
\end{aligned}
\end{equation*}
where the first inequality follows from the Cauchy-Schwarz inequality, which implies the known property
\begin{equation}\label{equ:inf_l2}
    \|p\|_{\infty} \leq C_n \|p\|_2 \qquad \forall p \in \mathbb{P}_n(\Omega).
\end{equation}

Since the quadrature rule has algebraic degree of precision $2n$, $\mathcal{L}_n f \in \mathbb{P}_{n}(\Omega)$ can be estimated in the uniform norm, that is,
\begin{equation*}
    \begin{aligned}
        \| \mathcal{L}_n f\|_{\infty} &\leq C_n \|\mathcal{L}_n f \|_2 = C_n \sqrt{\langle \mathcal{L}_n f, \mathcal{L}_n f \rangle} =  C_n \sqrt{\langle \mathcal{L}_n f, \mathcal{L}_n f \rangle_N}, \\
        & \leq C_n \sqrt{\langle  f, f \rangle_N} \leq C_n V^{1/2} \| f\|_{\infty}.
    \end{aligned}
\end{equation*}
For any $p \in \mathbb{P}_{n}(\Omega)$, we have $\mathcal{L}_n p =p$ and thus
\begin{equation*}
\| \mathcal{L}_n f - f  \|_{\infty} \leq \|  \mathcal{L}_n (f-p) \|_{\infty}  + \|   f-p \|_{\infty} \leq (1+C_n V^{1/2}) \|  f-p \|_{\infty}.
\end{equation*}
Let $p=p^{\ast}$ be the best approximation of $f$ in $\mathbb{P}_{n}(\Omega)$. Then we finish the proof. 
\end{proof}

Lemma \ref{lem:hyperL2UNI} provides a rough estimate for the hyperinterpolation operator $\mathcal{L}_n$. Some refined results  have been explored in various contexts: on the square \cite{caliari2007hyperinterpolation}, on the unit disc \cite{hansen2009norm}, on the sphere \cite{MR2274179}, on the spherical triangle \cite{sommariva2021spheritri} and in the cube \cite{caliari2008hyperinterpolation}. 
Additionally, estimates for the reciprocal of Christoffel functions $G_n(\mathbf{x},\mathbf{x})$ concerning the disc, ball, square, and cube are discussed in \cite{sommariva2014multivariate}.

The following result provides an upper bound on the uniform norm of the hard thresholding hyperinterpolation operator and estimates the uniform error associated with recovering the test function from noisy data using hard thresholding hyperinterpolation.

\begin{theorem}\label{thm:UpperBoundHardHyper_Uniform_Noise}
Adopt conditions of Theorem \ref{thm:ExactSolution}. Let $C_n= \sqrt{\max\limits_{\mathbf{x} \in \Omega} G_{n}(\mathbf{x},\mathbf{x})}$ with $G_n$ defined by \eqref{equ:Christoffel}. Then
\begin{equation}\label{equ:OperatorNormHardHyper_Uniform}
    \|\mathcal{H}_{n}^{\lambda} f\|_{\infty} \leq C_n V^{1/2}\|f\|_{\infty}
\end{equation}
and
\begin{equation}\label{equ:UpperBoundHardHyper_Uniform_Noise}
\|\mathcal{H}_n^{\lambda} f^{\epsilon} - f\|_{\infty}  \leq \|\mathcal{H}_n^{\lambda}f^{\epsilon} - \mathcal{L}_nf^{\epsilon} \|_{\infty} + (1+C_n)V^{1/2}E_{n}(f)+ C_nV^{1/2}\|{\epsilon} \|_{\infty}.
\end{equation}
\end{theorem}

\begin{proof}
By Lemma \ref{lem:hyperL2UNI} and Theorem \ref{thm:Noise}, we immediately obtain \eqref{equ:OperatorNormHardHyper_Uniform}. Next, we deduce
\begin{equation*}
    \|\mathcal{H}_n^{\lambda} f^{\epsilon} - f\|_{\infty} \leq \| \mathcal{H}_n^{\lambda} f^{\epsilon} - \mathcal{L}_nf^{\epsilon} \|_{\infty} + \|\mathcal{L}_nf^{\epsilon} -f \|_{\infty}.
\end{equation*}
By Lemma \ref{lem:hyperL2UNI} and the linearity of $\mathcal{L}_n$, we find that 
\begin{equation*}
\begin{aligned}
     \|\mathcal{L}_nf^{\epsilon} -f \|_{\infty}&= \|\mathcal{L}_n(f^{\epsilon}-f) + \mathcal{L}_nf -f\|_{\infty}, \\ 
     &\leq \|\mathcal{L}_n(f^{\epsilon}-f)\|_{\infty} + \| \mathcal{L}_nf -f\|_{\infty}, \\ 
     & \leq C_nV^{1/2}\|{\epsilon} \|_{\infty}+(1+C_n)V^{1/2}E_{n}(f) .
\end{aligned}
\end{equation*}
Therefore, we have completed the proof. 
\end{proof}

\subsection{Error analysis in practice}\label{sec:computational}
We further focus on the practical denoising capabilities of hard thresholding and Lasso hyperinterpolations, exploring how the choice of the regularization parameter $\lambda$ impacts their denoising effectiveness. Given $f \in \mathcal{C}(\Omega)$ and a specific noise $\epsilon$, we can collect sampling data $\mathbf{f}^{\epsilon}=[f^{\epsilon}(\mathbf{x}_1),\ldots,f^{\epsilon}(\mathbf{x}_N)]^{\rm{T}} \in \mathbb{R}^{N}$ from the contaminated function $f^{\epsilon}=f+\epsilon$.

\begin{theorem}\label{thm:priorparameter}
Adopt conditions of Theorem \ref{thm:ExactSolution}. Let all $\{\mu_{\ell}\}_{\ell=1}^{d_n}$ be 1 in \eqref{equ:lassohyperinterpolation}, and let $\bm{\epsilon}=[\epsilon_1,\ldots,\epsilon_N]^{\rm{T}} \in \mathbb{R}^{N}$, $\bm{\alpha}=[\alpha_1, \ldots, \alpha_{d_n}]^{\rm{T}} \in \mathbb{R}^{d_n}$ with $\alpha_{\ell} = \discreteinner{f^{\epsilon}, \Phi_{\ell}}$ and $\bm{\xi}=[\xi_1, \ldots, \xi_{d_n}]^{\rm{T}} \in \mathbb{R}^{d_n}$ with $\xi_{\ell} = \langle\epsilon, \Phi_{\ell} \rangle_N$ for $\ell=1, \ldots, d_n$. Define
\begin{equation*}
    J(\mathbf{z}):= \mathbf{z}^{\rm{T}}\mathbf{z}-2\mathbf{z}^{\rm{T}}\bm{\alpha}, \quad H(\mathbf{z}) := 2\mathbf{z}^{\rm{T}} \bm{\xi}, \quad \text{and} \quad E(\mathbf{z}) :=  \|\mathbf{W}^{1/2}(\mathbf{Az} - \mathbf{f}) \|_2^2
\end{equation*}
for $\mathbf{z}=[z_1, \ldots,z_{d_n}]^{\rm{T}} \in \mathbb{R}^{d_n}$. Then, the following results hold:
\begin{enumerate}
\item[(a)] For Lasso hyperinterpolation
\begin{equation*}
 E(\bm{\gamma}^{\lambda}) = J(\bm{\gamma}^{\lambda}) + H(\bm{\gamma}^{\lambda})+ \|\mathbf{W}^{1/2}\mathbf{f}\|_2^2,   
\end{equation*}
and the value of $J(\bm{\gamma}^{\lambda})$ is non-positive, where $\bm{\gamma}^{\lambda}=[\gamma_1^{\lambda},\ldots,\gamma_{d_n}^{\lambda}]^{\rm{T}} \in \mathbb{R}^{d_n}$ with $\gamma_{\ell}^{\lambda}=\eta_{S}(\alpha_{\ell},\lambda)$ for $\ell=1,\ldots,d_n$.
\item[(b)] For hard thresholding hyperinterpolation
\begin{equation*}
  E(\bm{\beta}^{\lambda}) = J(\bm{\beta}^{\lambda}) + H(\bm{\beta}^{\lambda})+ \|\mathbf{W}^{1/2}\mathbf{f}\|_2^2,  
\end{equation*}
and the value of $J(\bm{\beta}^{\lambda})$ is non-positive, where $\bm{\beta}^{\lambda}=[\beta_1^{\lambda},\ldots,\beta_{d_n}^{\lambda}]^{\rm{T}} \in \mathbb{R}^{d_n}$ with $\beta_{\ell}^{\lambda}=\eta_{H}(\alpha_{\ell},\lambda)$ for $\ell=1,\ldots,d_n$.
\item[(c)] If the regularization parameter $\lambda$ in Lasso and hard thresholding hyperinterpolations satisfies
\begin{equation}\label{cond:hardgreat}
       \sum_{\ell=1, |\alpha_{\ell}| > \lambda}^{d_n} \left[ \lambda  -2\rm{sgn}(\alpha_{\ell}) \xi_{\ell}  \right] \geq 0,
\end{equation}
then 
\begin{equation*}
    E(\bm{\beta}^{\lambda})   \leq E(\bm{\gamma}^{\lambda}).
\end{equation*}
\end{enumerate}
\end{theorem}

\begin{proof}
(a) For Lasso hyperinterpolation, since ${\mathbf{A}}^{\text{T}}{\mathbf{W}}{\mathbf{A}}={\mathbf{I}}_{d_n}$ is a $d_n \times d_n$ identity matrix, we easily have 
\begin{equation*}
\inner{\mathbf{W}^{1/2}\mathbf{A}\bm{\gamma}^{\lambda},\mathbf{W}^{1/2}\mathbf{A}\bm{\gamma}^{\lambda}}=\|\bm{\gamma}^{\lambda}\|_2^2.    
\end{equation*}
Next, since $\bm{\alpha}=\mathbf{A}^{\text{T}}\mathbf{W}\mathbf{f}^{\epsilon}$ and $\bm{\xi}=\mathbf{A}^{\text{T}}\mathbf{W}\bm{\epsilon}$, we can decompose $E(\bm{\gamma}^{\lambda})$ as the following
\begin{equation*} 
\begin{aligned}
E(\bm{\gamma}^{\lambda})
=&\|\bm{\gamma}^{\lambda}\|_2^2 - 2\inner{ \bm{\gamma}^{\lambda}, \mathbf{A}^{\text{T}}\mathbf{W}(\mathbf{f} + \bm{\epsilon} - \bm{\epsilon}) } + \|\mathbf{W}^{1/2}\mathbf{f}\|_2^2, \\
=& \|\bm{\gamma}^{\lambda}\|_2^2 - 2\inner{ \bm{\gamma}^{\lambda}, \mathbf{A}^{\text{T}}\mathbf{W}\mathbf{f}^{\epsilon} } + 2\inner{\bm{\gamma}^{\lambda}, \mathbf{A}^{\text{T}}\mathbf{W}\bm{\epsilon}} + \|\mathbf{W}^{1/2}\mathbf{f}\|_2^2, \\
=& \|\bm{\gamma}^{\lambda}\|_2^2 - 2\inner{ \bm{\gamma}^{\lambda}, \bm{\alpha}}   + 2\inner{\bm{\gamma}^{\lambda}, \mathbf{A}^{\text{T}}\mathbf{W}\bm{\epsilon}} + \|\mathbf{W}^{1/2}\mathbf{f}\|_2^2, \\
= & J(\bm{\gamma}^{\lambda}) +  H(\bm{\gamma}^{\lambda})+ \|\mathbf{W}^{1/2}\mathbf{f}\|_2^2 .
\end{aligned}
\end{equation*}
Since the coefficients $\{\gamma_{\ell}^{\lambda}\}_{\ell=1}^{d_n}$ of Lasso hyperinterpolation are processed by a soft thresholding operator, i.e.,
\begin{equation*}
    \gamma_{\ell}^{\lambda}= \eta_{S}(\alpha_{\ell}, \lambda) = \left\{  \begin{array}{cl} 
    \alpha_{\ell} - \lambda, & \text{ if } \alpha_{\ell}>\lambda, \\
    0, & \text{ if } |\alpha_{\ell}| \leq \lambda, \\
    \alpha_{\ell} + \lambda, & \text{ if } \alpha_{\ell}<\lambda, \\
    \end{array}\right.
\end{equation*}
we obtain
\begin{eqnarray*}
   J(\bm{\gamma}^{\lambda})= \summ{\ell}{d_n} \squarebrackets{\brackets{ \gamma^{\lambda}_{\ell}}^2 -2\gamma^{\lambda}_{\ell}\alpha_{\ell} } \leq 0.
\end{eqnarray*}

(b) Taking the same techniques as in (a) gives 
\begin{equation*}
     E(\bm{\beta}^{\lambda}) = J(\bm{\beta}^{\lambda}) + H(\bm{\beta}^{\lambda})+ \|\mathbf{W}^{1/2}\mathbf{f}\|_2^2.
\end{equation*}
Next, according to the fact that all coefficients $\{\beta_{\ell}^{\lambda}\}_{\ell=1}^{d_n}$ of hard thresholding hyperinterpolation are given by 
\begin{equation*}
        \beta^{\lambda}_{\ell}=\eta_H(\alpha_{\ell},\lambda) = \left\{\begin{array}{cl}
        \alpha_{\ell}, & \text{ if } |\alpha_{\ell}| > \lambda, \\
        0 , & \text{ if }  |\alpha_{\ell}| \leq \lambda,
    \end{array}  \right.
\end{equation*}
we have $(\beta^{\lambda}_{\ell})^2 -2\beta^{\lambda}_{\ell} \alpha_{\ell} \leq 0$ for $\ell=1,\ldots, d_n$.

(c) %Let $\mathds{1}_{|\alpha_{\ell}|>\lambda}(\alpha_{\ell})$ be an indicator function defined as
%\begin{equation*}
%    \mathds{1}_{|\alpha_{\ell}|>\lambda} (\alpha_{\ell}):= \left\{\begin{array}{cl} 
%        1, & \text{ if } |\alpha_{\ell}|>\lambda, \\
%        0, & \text{ if } |\alpha_{\ell}|\leq \lambda.
%    \end{array}\right.
%\end{equation*}
Inspired by the proof of (a) and (b), we deduce that
\begin{equation*}\label{equ:ErrHard}
    E(\bm{\beta}^{\lambda}) = \sum\limits_{\ell=1, |\alpha_{\ell}| > \lambda }^{d_n} (2\alpha_{\ell}\xi_{\ell} - \alpha_{\ell}^2 ) + \sum\limits_{j=1}^{N}w_jf^{2}(\mathbf{x}_j),
\end{equation*}
and
\begin{equation}\label{equ:ErrLasso}
\begin{aligned}
E(\bm{\gamma}^{\lambda}) =& \sum\limits_{\ell=1, |\alpha_{\ell}| > \lambda }^{d_n} (2\alpha_{\ell}\xi_{\ell} - \alpha_{\ell}^2 ) + \sum\limits_{j=1}^{N}w_jf^{2}(\mathbf{x}_j)  \\
                               & + \sum\limits_{\ell=1, |\alpha_{\ell}| > \lambda }^{d_n} [\lambda^2  - 2 \rm{sgn}( \alpha_{\ell})\lambda\xi_{\ell}].                 
\end{aligned}                               
\end{equation}
Since $\lambda$ is positive, then $E(\bm{\beta}^{\lambda}) \leq E(\bm{\gamma}^{\lambda})$.  
\end{proof}

\begin{remark}
For hyperinterpolation $\mathcal{L}_n f^{\epsilon}=\sum_{\ell=1}^{d_n}\alpha_{\ell}\Phi_{\ell}=\sum_{\ell=1}^{d_n}\langle f^{\epsilon}, \Phi_{\ell}\rangle_N\Phi_{\ell}$, we have 
\begin{equation*}
        E(\bm{\alpha}) = \sum\limits_{\ell=1 }^{d_n} (2\alpha_{\ell}\xi_{\ell} - \alpha_{\ell}^2 ) + \sum\limits_{j=1}^{N}w_jf^{2}(\mathbf{x}_j).
\end{equation*}
Observing equation \eqref{equ:ErrLasso}, one can find that Lasso hyperinterpolation not only selects some basis but also introduces an additional term 
\begin{equation*}
    \sum\limits_{\ell=1, |\alpha_{\ell}| > \lambda }^{d_n} [\lambda^2  - 2 \rm{sgn}( \alpha_{\ell})\lambda\xi_{\ell}].
\end{equation*}
We refer to $\{\xi_{\ell}\}_{\ell=1}^{d_n}$ \emph{noise coefficients}, which are derived from hyperinterpolation approximating the noise $\epsilon$.
\end{remark}

\section{Examples}\label{sec:examples}
In this section we test the classical hyperinterpolation as well as Lasso and hard thresholding hyperinterpolations over  the sphere ${\mathbb{S}}^2=\{{\mathbf{x}}: \|{\mathbf{x}}\|_2=1\}$ , the spherical triangle $\mathcal{T}$ in $\mathbb{S}^2 \subset \mathbb{R}^3$ and the cube. Differently from the experiment on the cube and sphere, an orthonormal basis is not theoretically available and must be computed numerically over the spherical triangle.

Depending on the numerical experiments, we have added independent noise to the evaluation of a function $f$ on the $N$ nodes $\{ {\mathbf{x}_j}\}_{j=1}^{N}$. In particular, we considered
\begin{itemize}
\item[$\bullet$] {\em{Gaussian noise}} ${\cal{N}}(0,\sigma^2)$ from a normal distribution with mean 0
        and standard deviation {\tt{sigma}}=$\sigma$, implemented via the Matlab command
\begin{center}
{\tt{sigma*randn(N,1)}}.
\end{center}
\item[$\bullet$]  {\em{Impulse noise}} ${\cal{I}}(a)$ that takes a uniformly distributed random values in $[-a,a]$ with probability density $1/(2a)$ by means of the Matlab command
\begin{center}
{\tt{a*(1-2*rand(N,1)).*binornd(1,0.5,N,1)}},
\end{center}
where {\tt{binornd(1,0.5,N,1)}} generates an array of $N \times 1$ random binary numbers (0 or 1), with each number having the probability $1/2$ of being 1 and the probability $1/2$ of being 0.
\end{itemize}

The regularization parameters $\lambda$ for hard thresholding hyperinterpolation $\mathcal{H}_n^{\lambda} f^{\epsilon}$ and Lasso hyperinterpolation $\mathcal{L}_n^{\lambda} f^{\epsilon}$ are computed for a sequence of increasing smoothing parameters, specifically $\lambda=2^{s}$ where $s$ ranges from $-15$ to $7$ in increments of $0.1$. The optimal $\lambda$ values for both $\mathcal{H}_n^{\lambda} f^{\epsilon}$ and  $\mathcal{L}_n^{\lambda} f^{\epsilon}$ are selected to minimize the $L_2$ errors, respectively. More examples, including the interval and the polygon, can be found on the website:\\
\href{https://github.com/JiaShuRan/HardThresholdingHyperinterpolation}{https://github.com/JiaShuRan/HardThresholdingHyperinterpolation}.

\subsection{The sphere}\label{sec:sphere}
As the first domain we consider the unit-sphere $\Omega$, i.e., $\Omega:={\mathbb{S}}^2=\{{\mathbf{x}}: \|{\mathbf{x}}\|_2=1\}$  with the surface area  $V=\int_{\mathbb{S}^2} \omega(\mathbf{x}){\rm{d}}\mathbf{x}=4\pi$. Let $\mathbb{P}_n(\mathbb{S}^2)$ be the space of spherical polynomials of degree at most $n$ with dimension $d_n = \dim \mathbb{P}_n (\mathbb{S}^2)= (n+1)^2$. We adopt the so-called {\em{spherical harmonics}} as the orthonormal polynomial basis {\cite{atkinson2012spherical}}. In 1977, Delsarte, Goethals and Seidel introduced the the famous {\em{spherical $t$-designs}} on ${\mathbb{S}}^2$ in the pioneering work {\cite{delsarte1977spherical}}.  A point set $\{{\mathbf{x}}_1,\ldots,{\mathbf{x}}_N\} \subset {\mathbb{S}}^2$ is a spherical $t$-design if it satisfies
\begin{equation*} 
    {\frac{1}{N}} \sum_{j=1}^N p({\mathbf{x}}_j) = \frac{1}{4\pi} \int_{{\mathbb{S}}^2} p({\mathbf{x}}) \text{d} {\omega}({\mathbf{x}}) \qquad \forall p \in {\mathbb{P}}_t({\mathbb{S}}^2).
\end{equation*}
In the following, we shall consider in particular an {\em{efficient spherical design}} proposed in {\cite{womersley2018efficient}} with a relatively small number of quadrature points $N$ and a uniformly bounded ratio of the covering radius to the packing radius. In our tests, since we intend to compute several type of hyperinterpolants of total degree at most $n=15$, it is necessary to adopt a rule with algebraic degree of exactness $2n=30$, consisting of $N=482$ points.
We thus have that for any $f,g \in {\mathbb{P}}_n({\mathbb{S}}^2)$
\begin{equation*}
    \int_{{\mathbb{S}}^2} f({\mathbf{x}})g({\mathbf{x}}) \text{d} {\omega}({\mathbf{x}})=  \langle f,g \rangle=\langle f,g \rangle_{482}:=\sum_{j=1}^{482} w_j f({\mathbf{x}}_j) g({\mathbf{x}}_j),
\end{equation*}
where the positive quadrature weights $w_1 = \cdots = w_{482}= 4\pi /482$.

Following {\cite{an2021lasso}} and letting $\mathbf{z}_1=[1,0,0]^{\rm{T}}$, $\mathbf{z}_2=[-1,0,0]^{\rm{T}}$, $\mathbf{z}_3=[0,1,0]^{\rm{T}}$, $\mathbf{z}_4=[0,-1,0]^{\rm{T}}$, $\mathbf{z}_5=[0,0,1]^{\rm{T}}$, and $\mathbf{z}_6=[0,0,-1]^{\rm{T}}$, we considered the function
\begin{equation*}
    f({\mathbf{x}})=  {\frac{1}{3}} \sum_{i=1}^6 \Phi_2(\| {\mathbf{z}}_i - {\mathbf{x}}\|_2),
\end{equation*}
where $\Phi_2(r):={\tilde{\Phi}}_2(r/\delta_2)$ and $\| {\mathbf{z}}_i - {\mathbf{x}}\|_2 = [(\mathbf{z}_i-\mathbf{x})^{\rm{T}}(\mathbf{z}_i-\mathbf{x})]^{\frac{1}{2}}$, in which ${\tilde{\Phi}}_2$ is the $\mathcal{C}^6$ compactly supported of minimal degree Wendland function defined as
\begin{equation*}
    {\tilde{\Phi}}_2(r):=(\max\{1-r,0\})^6 (35 r^2 + 18r +3)
\end{equation*}
and $\delta_2=\frac{9 \Gamma(5/2)}{2\Gamma(3)} $. 

\begin{figure}[htbp]
    \centering
    \includegraphics[width=\textwidth]{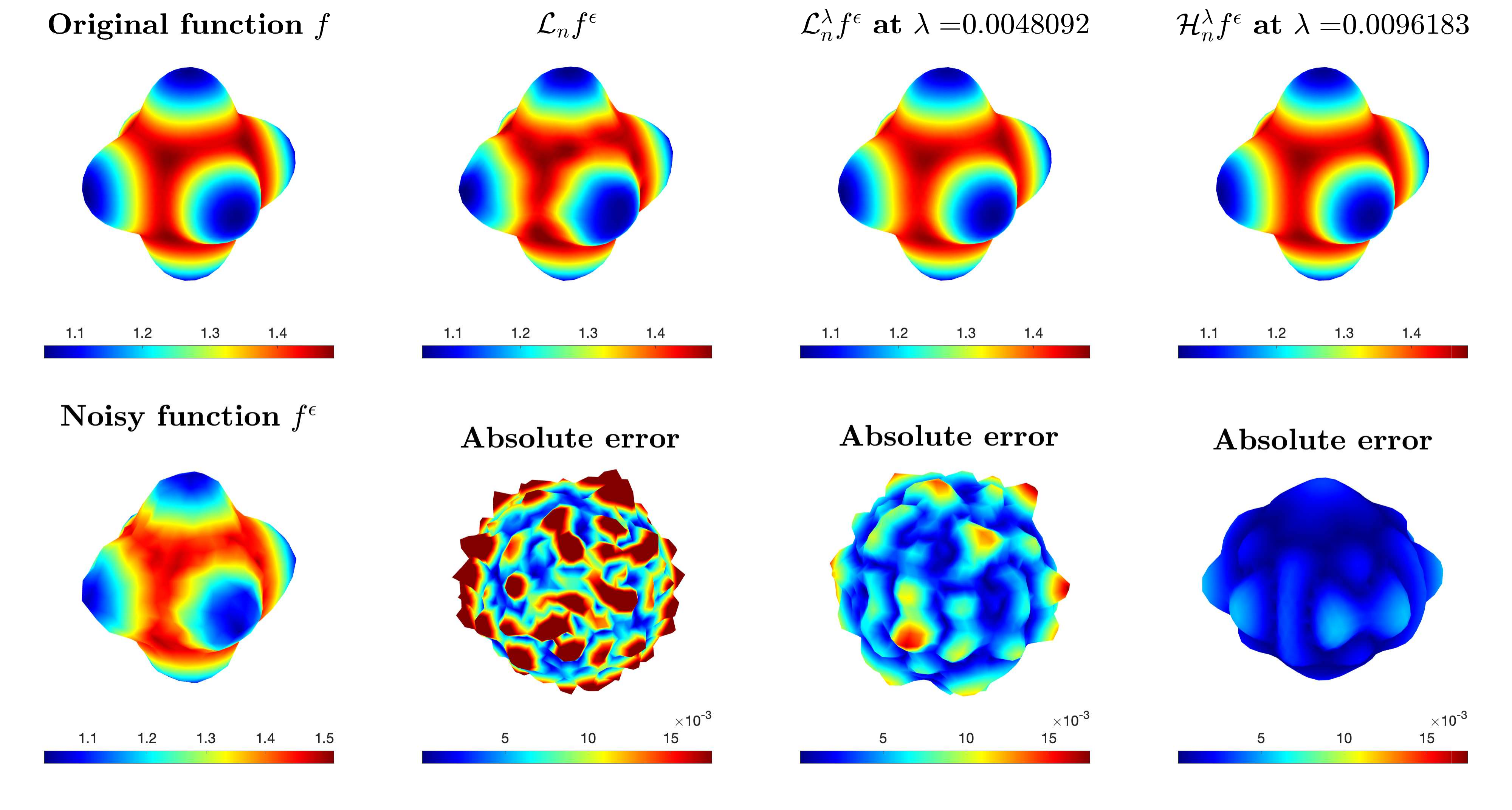}
    \caption{Approximate  $f(\mathbf{x}) ={\frac{1}{3}} \sum_{i=1}^6 \Phi_2(\| {\mathbf{z}}_i - {\mathbf{x}}\|_2)$, perturbed by impulse noise ($a=0.02$) and Gaussian noise ($\sigma=0.02$), over the unit-sphere ${\mathbb{S}}^2$, via
    hyperinterpolation $\mathcal{L}_n f^{\epsilon}$,
    Lasso hyperinterpolation $\mathcal{L}_{n}^{\lambda}f^{\epsilon}$, and
    hard thresholding hyperinterpolation $\mathcal{H}_{n}^{\lambda} f^{\epsilon}$ at $n=15$.}
    \label{sphere_denoise}
\end{figure}

\begin{figure}[htbp]
    \centering
    \includegraphics[width=\textwidth]{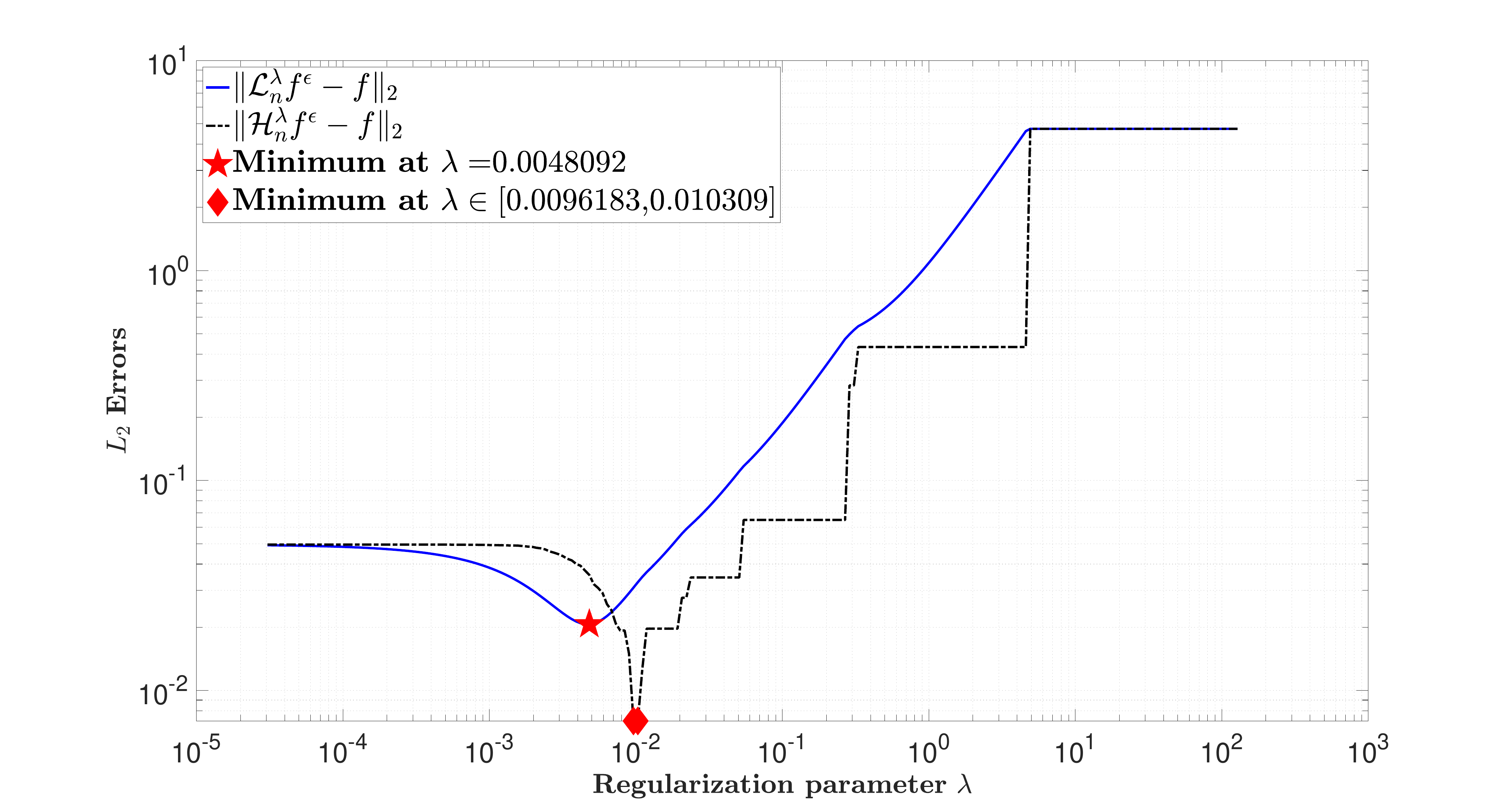}
    \caption{The choices of regularization parameter $\lambda$  for Lasso hyperinterpolation $\mathcal{L}_n^{\lambda} f^{\epsilon}$ and hard thresholding hyperinterpolation   $\mathcal{H}_n^{\lambda} f^{\epsilon}$ at $n=15$ approximating
    $f(\mathbf{x}) ={\frac{1}{3}} \sum_{i=1}^6 \Phi_2(\| {\mathbf{z}}_i - {\mathbf{x}}\|_2)$, perturbed by impulse noise ($a=0.02$) and Gaussian noise ($\sigma=0.02$), over the unit-sphere ${\mathbb{S}}^2$.}
    \label{sphere}
\end{figure}

In Figures \ref{sphere_denoise} and \ref{sphere}, we fixed the quadrature exactness $2n=30$, quadrature points $N=482$, Gaussian noise ${\cal{N}}(0,\sigma^2)$ with standard deviation $\sigma=0.02$ and impulse noise ${\cal{I}}(a)$ relatively to the level $a=0.02$.  Figure \ref{sphere_denoise} compares the denoising effectiveness of three methods: hard thresholding hyperinterpolation $\mathcal{H}_n^{\lambda} f^{\epsilon}$, Lasso hyperinterpolation $\mathcal{L}_n^{\lambda} f^{\epsilon}$ and classical hyperinterpolation $\mathcal{L}_n f^{\epsilon}$, all at degree $n=15$. Both hard thresholding and Lasso hyperinterpolation outperform classical hyperinterpolation in reconstructing the test function. Specifically, as illustrated in Figure  \ref{sphere},  $\mathcal{H}_n^{\lambda} f^{\epsilon}$ achieves the lowest $L_2$ error of 0.00715667 at $\lambda\in [0.0096183,0.010309]$ with only 8 non-zero coefficients, while $\mathcal{L}_n^{\lambda} f^{\epsilon}$ reaches a minimum $L_{2}$ error of 0.02066131 at $\lambda=0.0048092$ but with 38 non-zero coefficients. 

Figure \ref{sphere} illustrates the $L_2$ errors as functions of the regularization parameter $\lambda$ for two methods: Lasso hyperinterpolation $\mathcal{L}_n^{\lambda} f^{\epsilon}$ (shown in blue) and hard thresholding hyperinterpolation $\mathcal{H}_n^{\lambda} f^{\epsilon}$ (shown in dashed black).
The dashed black curve, which has a staircase-like shape, is a direct result of the hard thresholding operator $\eta_H$ used in $\mathcal{H}_n^{\lambda} f^{\epsilon}$.  This operator retains only those coefficients whose absolute values exceed the regularization parameter $\lambda$ and sets all others to zero. This process leads to abrupt changes in the set of retained coefficients as $\lambda$ is adjusted, causing the error to increase sharply at certain points—hence, the staircase pattern.

In contrast, the blue curve representing the $L_2$ errors for $\mathcal{L}_n^{\lambda} f^{\epsilon}$ does not exhibit such abrupt changes. Instead, the Lasso method smoothly penalizes the coefficients proportionally, leading to a gradual and continuous variation in the error as $\lambda$ changes. This results in a much smoother curve compared to the staircase-like pattern observed with the hard thresholding approach.

In addition, we verify that the section where $\|\mathcal{H}_n^{\lambda} f^{\epsilon}-f\|_2 \leq \|\mathcal{L}_n^{\lambda} f^{\epsilon}-f\|_2$ in Figure \ref{sphere}  implies that the regularization parameter $\lambda$ satisfies condition \eqref{cond:hardgreat} stated in Theorem \ref{thm:priorparameter}.
%\begin{figure}[htbp]
%  \centering
%  \includegraphics[width=\textwidth]{figure2-eps-converted-to.pdf}\\
%  \caption{Approximation results of $ f({\mathbf{x}})=  {\frac{1}{3}} \sum_{i=1}^6 \Phi_2(\| {\mathbf{z}}_i - {\mathbf{x}}\|_2),$ perturbed by Gaussian noise with $\sigma = 0.02$ and single impulsive noise with $a=0.02$ over the sphere via hard thresholding hyperinterpolation $\mathcal{H}_{15}^{\lambda}{f^{\epsilon}}$ and Lasso hyperinterpolation $\mathcal{L}_{15}^{\lambda}{f^{\epsilon}}$, hyperinterpolation $\mathcal{L}_{15} f^{\epsilon}$  with the regularization parameter $\lambda^{\ast} = \lambda(9) = 0.0097$.}\label{sphere_denoise}
%\end{figure}

%\begin{figure}[htbp]
%  \centering
%  \includegraphics[width=\textwidth]{figure3-eps-converted-to.pdf}\\
%  \caption{Approximation results of $ f({\mathbf{x}})=  {\frac{1}{3}} \sum_{i=1}^6 \Phi_2(\| {\mathbf{z}}_i - {\mathbf{x}}\|_2),$ perturbed by Gaussian noise with $\sigma = 0.02$ and single impulsive noise with $a=0.02$ over the sphere  via hard thresholding hyperinterpolation $\mathcal{H}_{8}^{\lambda}{f^{\epsilon}}$ and Lasso hyperinterpolation $\mathcal{L}_{8}^{\lambda}{f^{\epsilon}}$ with quadrature points $N=482$. (a): the condition \eqref{cond:hardgreat} in Theorem \ref{thm:priorparameter}; (b): the relation between sparsity and $L_2$ errors.}\label{sphere}
%\end{figure}

\subsection{The spherical triangle}\label{sec:spt} 
\begin{figure}[htbp]
  \centering
  \includegraphics[scale=0.25,clip]{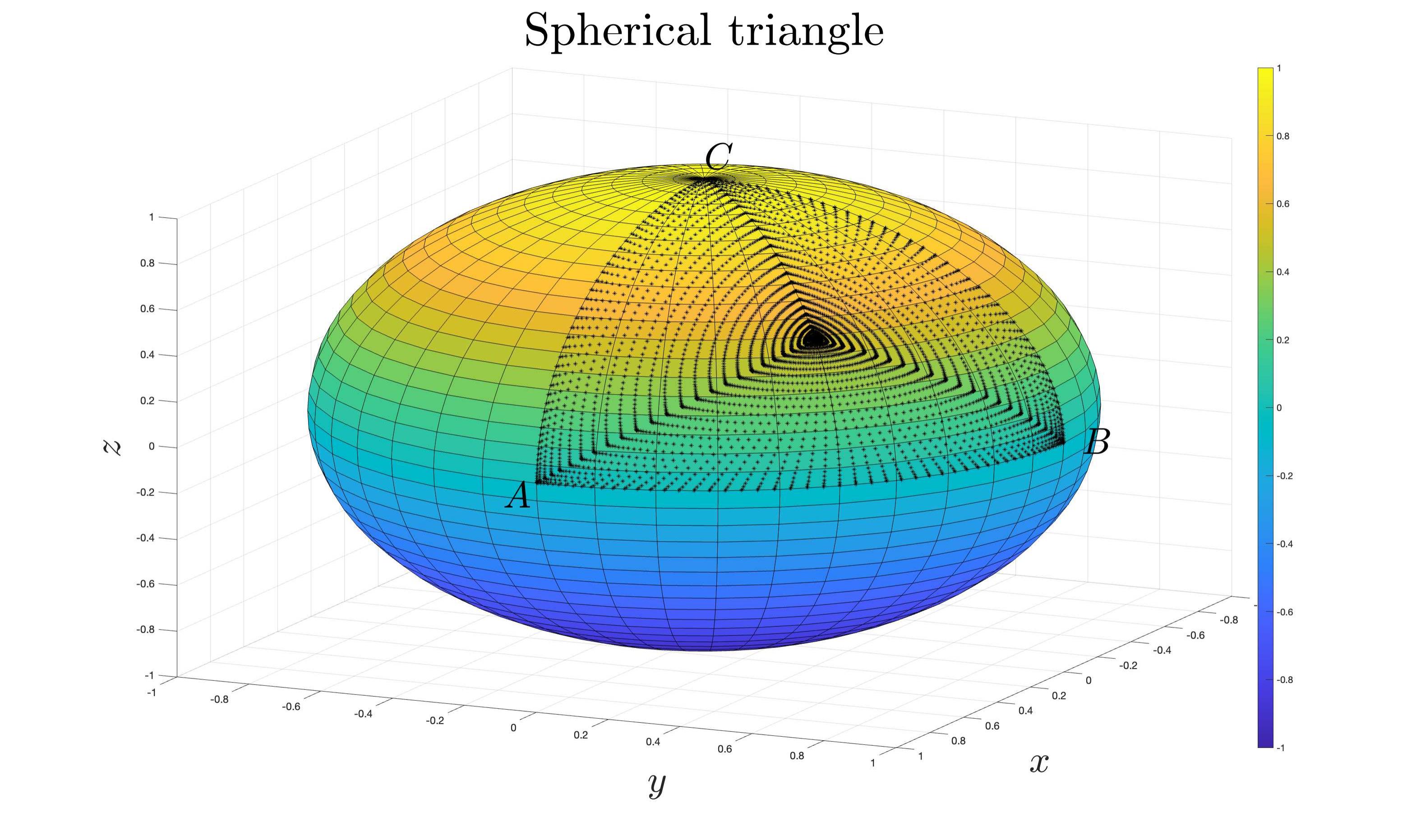}\\
  \caption{The domain $\mathcal{T}=\overset{\frown}{ABC}$ with vertices $A=[1,0,0]^{{\rm{T}}}$, $B=[0,1,0]^{{\rm{T}}},C=[0,0,1]^{{\rm{T}}}$ in which we perform our tests. We represent in black the $3N_{2n+m}=4620$ nodes of the cubature rule with $n=10$ and $m=34$.}\label{ST}
\end{figure}

In \cite{Sommariva2021Tchakaloff,sommariva2021spheritri}, authors pay attention to the numerical construction of hyperinterpolation over the spherical triangle, which plays an important role in geomathematics.  
Here, we apply hard thresholding hyperinterpolation to the spherical triangle $\mathcal{T}=\overset{\frown}{ABC}$ , namely the octant with vertices $A=[1,0,0]^{{\rm{T}}}$, $B=[0,1,0]^{{\rm{T}}},C=[0,0,1]^{{\rm{T}}}$ in Figure \ref{ST}. This spherical triangle $\mathcal{T}$ lies on the unit sphere $ \mathbb{S}^2\subset\mathbb{R}^3$, and its surface area can be calculated as follows: 
\begin{equation*}
   V=\int_{\mathcal{T}}\omega(\mathbf{x})\text{d}\mathbf{x} = \frac{1}{8} \int_{\mathbb{S}^2}\omega(\mathbf{x})\text{d}\mathbf{x} = \frac{\pi}{2},
\end{equation*}
where $\text{d}\omega=\omega(\mathbf{x})\text{d}\mathbf{x}$ and $\omega(\mathbf{x})$ is an area measure on the sphere $\mathbb{S}^2$.

Let $\mathbb{P}_n(\mathcal{T})$ be the polynomial space of 3-variate total-degree product Chebyshev basis $\{\phi_1,\ldots,\phi_{d_n}\}$  of a cartesian box containing the spherical triangle. The dimension of $\mathbb{P}_n(\mathcal{T})$ is $d_n = \dim(\mathbb{P}_n(\mathcal{T}))=(n+1)^2$. Following \cite{sloan1995hyperinterpolation}, in order to orthonormalize the basis with respect to the discrete measure generated by the quadrature formula, which by polynomial exactness up to degree 2n is orthonormal also with respect to ${\rm{d}}\omega$, we can compute the QR factorization with $\mathbf{Q} \in \mathbb{R}^{N \times d_n}$ and $\mathbf{R} \in\mathbb{R}^{d_n \times d_n}$, and construct the orthonormal basis $\{\Phi_1,\ldots,\Phi_{d_n}\}$ as
\begin{equation*}
    \mathbf{W}^{1/2}\mathbf{A} = \mathbf{QR}, \quad (\Phi_1,\ldots,\Phi_{d_n}) = (\phi_1,\ldots,\phi_{d_n})\mathbf{R}^{-1},
\end{equation*}
where $\mathbf{A} \in \mathbb{R}^{N \times d_n}$ with $[\mathbf{A}]_{j\ell} = \phi_{\ell}(\mathbf{x}_j)$ and $\mathbf{W}=\text{diag}(w_1,\ldots,w_N)$. Then we obtain the new Vandermonde-like matrix on orthonormal basis $\mathbf{U}=\mathbf{AR}^{-1} \in \mathbb{R}^{N \times d_n}$ with $[\mathbf{U} ]_{j \ell}= \Phi_{\ell}(\mathbf{x}_j)$.

\begin{figure}[htbp]
  \centering
  \includegraphics[scale=0.3,clip]{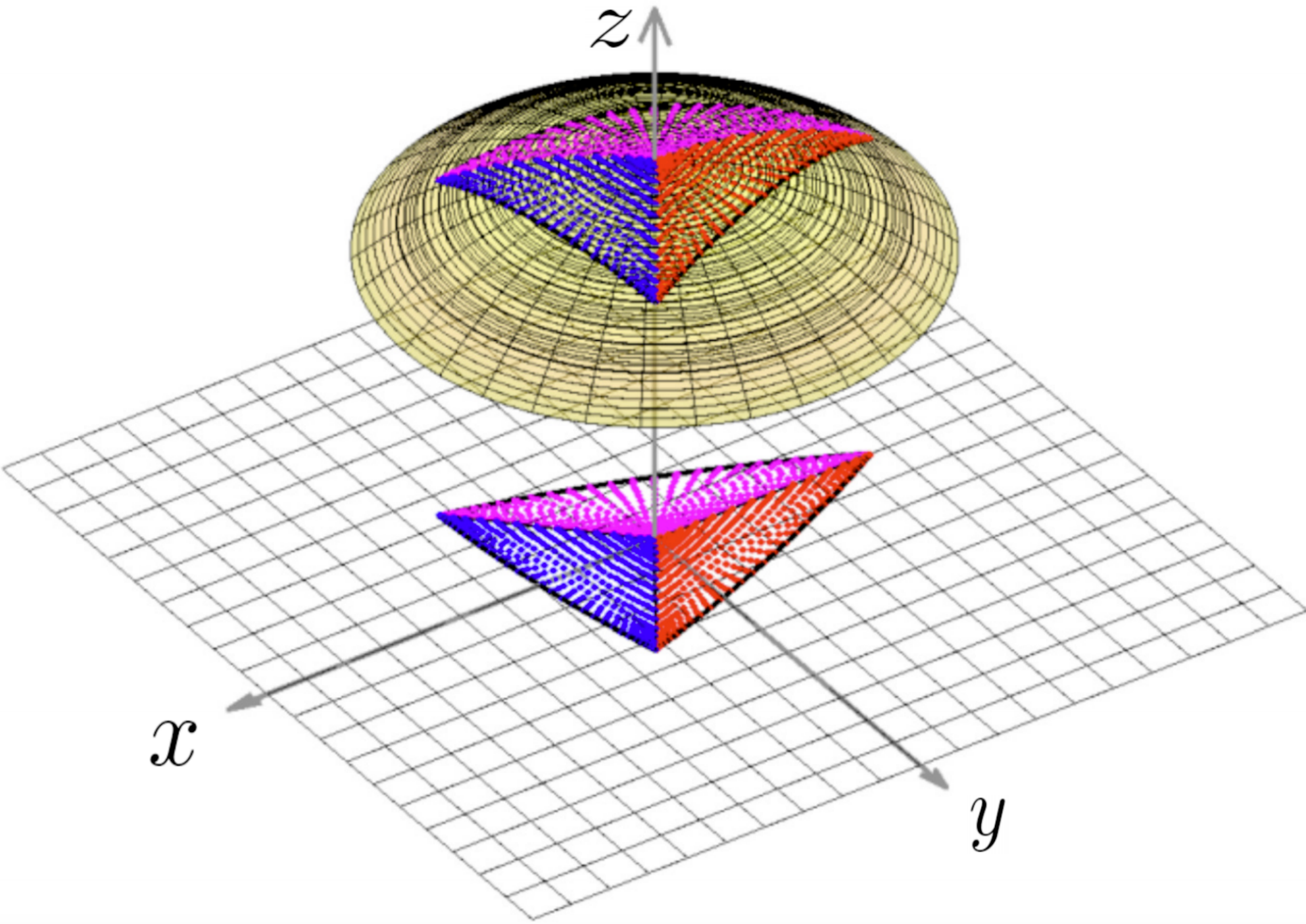}\\
  \caption{Quadrature nodes on a spherical triangle rotated with centroid at the north pole, completely contained in a hemisphere and lifted from the projected elliptical triangle, before compression.}\label{ST_split}
\end{figure}

There exists an algorithm for the computation of nodes and positive weights of a quadrature formula on spherical triangles, which is nearly exact for algebraic polynomials of a given degree $n$ on $\mathbb{S}^2$, and whose cardinality does not exceed the dimension of the corresponding polynomial space \cite{Sommariva2021Tchakaloff}. However, we only consider the general case, i.e., not compressing the quadrature nodes to a low-cardinality. To find a quadrature formula with positive weights and exactness degree $2n$ on $\mathcal{T}$, we can consider a general case, where the spherical triangle with centroid $(A+B+C)/\| A+B+C\|_2$ at the north pole. Then the spherical triangle can be projected on the equatorial plane as an ``elliptical triangle'' which can be split into three elliptical sectors, say $S_1,S_2,S_3$; see Figure \ref{ST_split}. For a continuous function $f(x,y,z)$ on $\mathcal{T}$, we have

\begin{equation*}
    \int_{\mathcal{T}} f(x,y,z) \text{d} \omega \approx \sum_{j=1}^{3N_{2n+m}} \frac{w_j}{\sqrt{1-x_j^2 - y_j^2}} f\left( x_j,y_j, \sqrt{1-x_j^2 - y_j^2}\right),
\end{equation*}
which is nearly exact for spherical polynomials of degree not exceeding $n$, and where $N_{2n+m}=(2n+m+1)\lceil (2n+m+2)/2 \rceil$ is the number of quadrature nodes $\{(x_j,y_j)\}$ on each elliptical sector, $\lceil \cdot \rceil$ is the ceiling function and $\{w_j\}$ are the corresponding positive weights, and $m$ denotes the degree of a suitable bivariate polynomial that approximates $1/\sqrt{1-x^2-y^2}$ at machine precision on the elliptical triangle which is the projection of $\mathcal{T}$ onto the $xy$-plane.

As the numerical experiment, we examine the reconstruction of the function 
\begin{equation*}
    f(x,y,z) = \exp\left( -\left(x -\frac{1}{\sqrt{3}}\right)^2- \left(y -\frac{1}{\sqrt{3}}\right)^2 - \left(z -\frac{1}{\sqrt{3}}\right)^2 \right),
\end{equation*}
to which we have added Gaussian noise ${\cal{N}}(0,\sigma^2)$ with standard deviation $\sigma=0.2$ and impulse noise ${\cal{I}}(a)$ relatively to the level $a=0.1$. 

\begin{figure}[htbp]
  \centering
  \includegraphics[width=\textwidth]{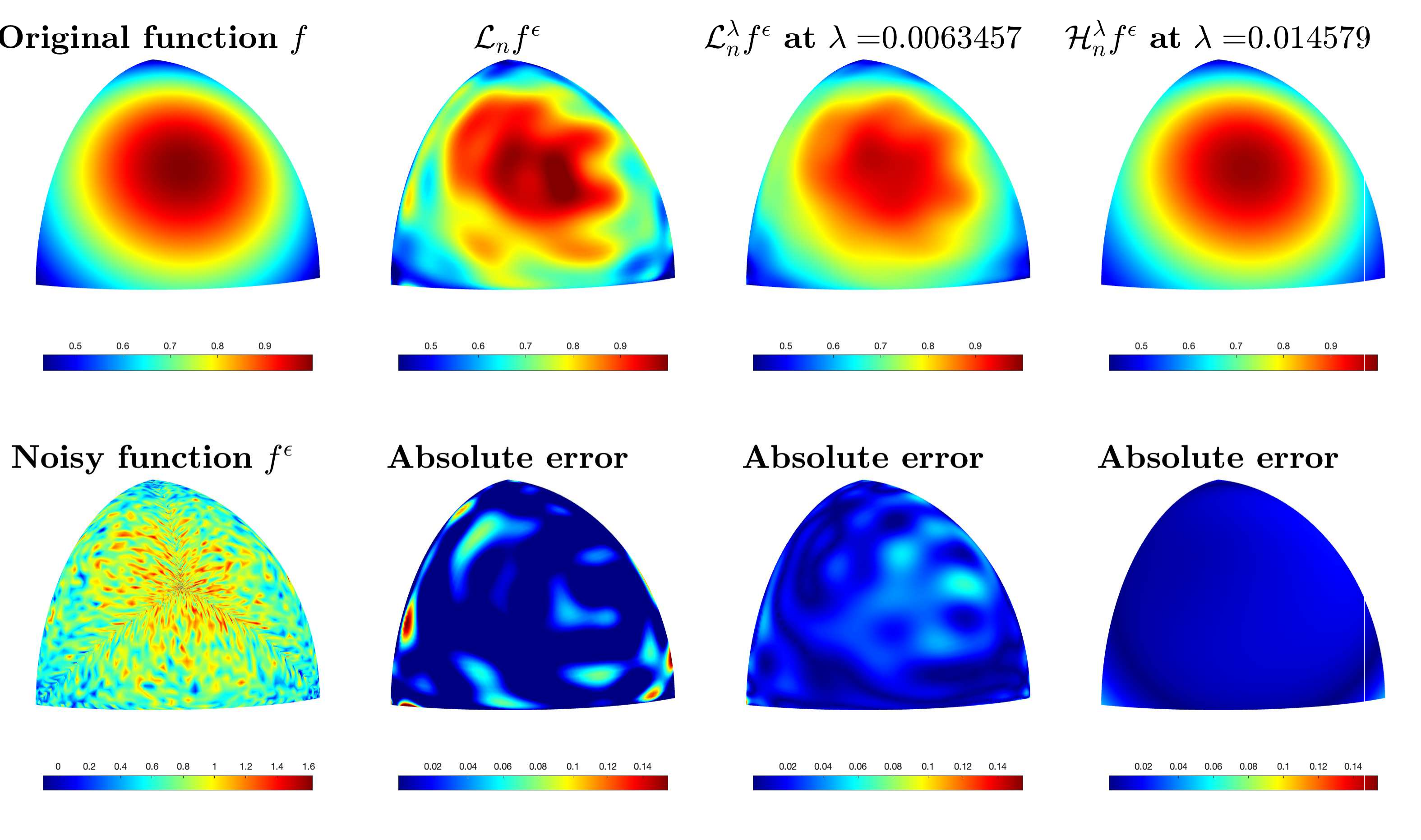}\\
  \caption{Approximate $f(x,y,z)=\text{exp}(-(x- 1/\sqrt{3})^2-(y- 1/\sqrt{3})^2-(z- 1/\sqrt{3})^2)$ perturbed by impulse noise ($a = 0.1$) and Gaussian noise ($\sigma = 0.2$) over $\mathcal{T}$ depicted in Figure {\ref{ST}} via hyperinterpolation $\mathcal{L}_{n}f^{\epsilon}$, Lasso hyperinterpolation $\mathcal{L}_{n}^{\lambda}{f^{\epsilon}}$ and hard thresholding hyperinterpolation $\mathcal{H}_{n}^{\lambda}{f^{\epsilon}}$ at $n=10$.}\label{spt_denoise}
\end{figure}

\begin{figure}[htbp]
  \centering
  \includegraphics[width=\textwidth]{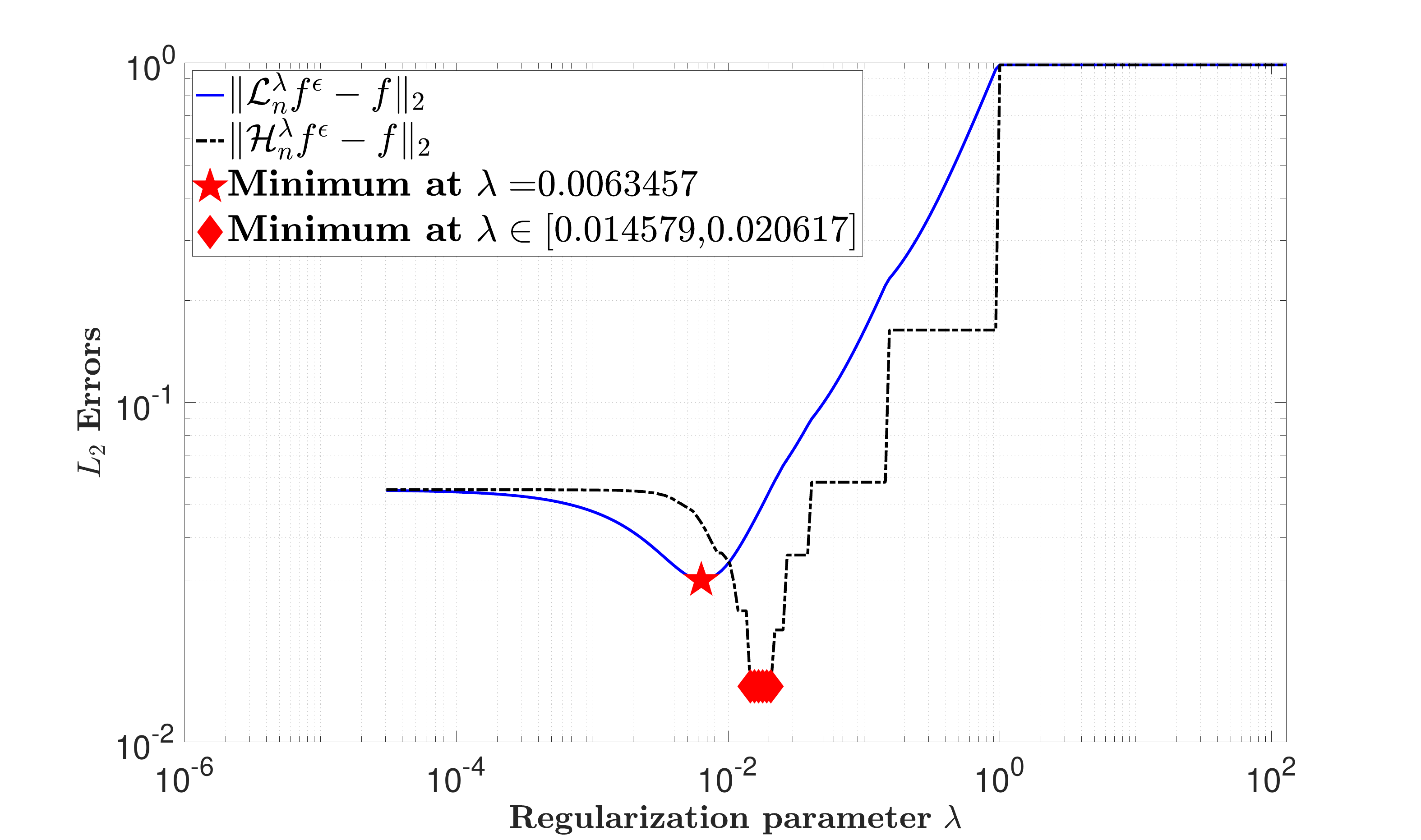}\\
  \caption{The choices of regularization parameter $\lambda$ for Lasso hyperinterpolation $\mathcal{L}_{n}^{\lambda}{f^{\epsilon}}$ and hard thresholding hyperinterpolation $\mathcal{H}_{n}^{\lambda}{f^{\epsilon}}$ at $n=10$ approximating $f(x,y,z)=\text{exp}(-(x- 1/\sqrt{3})^2-(y- 1/\sqrt{3})^2-(z- 1/\sqrt{3})^2)$, perturbed by impulse noise ($a = 0.1$) and Gaussian noise ($\sigma = 0.2$), over $\mathcal{T}$ depicted in Figure {\ref{ST}}.}\label{SphericalTriangle}
\end{figure}

Build on the analysis in section \ref{sec:sphere}, Figure \ref{spt_denoise} presents the comparison of denoising performance among hard thresholding hyperinterpolation $\mathcal{H}_n^{\lambda} f^{\epsilon}$, Lasso hyperinterpolation $\mathcal{L}_n^{\lambda} f^{\epsilon}$, and classical hyperinterpolation $\mathcal{L}_n f^{\epsilon}$ for $n=10$. Consistent with earlier results, $\mathcal{H}_n^{\lambda} f^{\epsilon}$ outperforms its counterparts, achieving the lowest $L_2$ errors with fewer non-zero coefficients. Specific details are depicted in Figure \ref{SphericalTriangle}, where $\mathcal{H}_n^{\lambda} f^{\epsilon}$ reaches an $L_2$ error of 0.014592 is obtained for $\lambda\in [0.014579,0.020617]$ with 5 non-zero coefficients, while $\mathcal{L}_n^{\lambda}f^{\epsilon}$ achieves its minimum $L_{2}$ error of 0.029886 at $\lambda=0.0063457$ with 26 non-zero coefficients.

%In Figures \ref{spt_denoise} and \ref{SphericalTriangle}, we fixed the quadrature exactness $2n=20$, quadrature points $N=3N_{2n+m}=4620$.  In Figure \ref{spt_denoise}, we test the denoising abilities of hard thresholding hyperinterpolation $\mathcal{H}_n^{\lambda} f^{\epsilon}$, Lasso hyperinterpolation $\mathcal{L}_n^{\lambda} f^{\epsilon}$ with the given regularization parameter $\lambda$ to reach its smallest $L_{2}$ error, and classical hyperinterpolation $\mathcal{L}_n f^{\epsilon}$, where the degree is $10$.  As illustrated in Figure \ref{SphericalTriangle} for $\mathcal{H}_n^{\lambda} f^{\epsilon}$, the smallest $L_2$ error of 0.004817 is obtained for $\lambda\in [0.011842,0.017948]$ with only 5 non-zero coefficients, while the smallest $L_{2}$ error of 0.017053 for $\mathcal{L}_n^{\lambda}f^{\epsilon}$ is obtained for $\lambda=0.00634572$ with 27 non-zero coefficients.

\subsection{The cube}\label{sec:cube} 
As the third domain we consider the cube $\Omega=[-1,1]^3 \subset {\mathbb{R}}^3$.  For two functions $f,g \in \mathcal{C}(\Omega)$ with respect to a certain measure such that $V=\int_{[-1,1]^3}  w({\mathbf{x}}) \text{d}{\mathbf{x}}=1$, we can define a scalar product
\begin{eqnarray*}
    \langle f,g \rangle = \int_{[-1,1]^3} f({\mathbf{x}}) g({\mathbf{x}}) w({\mathbf{x}}) \text{d}{\mathbf{x}}, {\hspace{0.5cm}}
w({\mathbf{x}}) = {\frac{1}{\pi^3}} \prod {\frac{1}{\sqrt{1-x_i^2}}},
\end{eqnarray*}
where as usual ${\mathbf{x}}=[x_1,x_2,x_3]^{\rm{T}} \in \mathbb{R}^3$.

Let $\mathbb{P}_n(\Omega)$ be the polynomial space of the product Chebyshev orthonormal basis $\{\Phi_{\ell}\}_{\ell=1}^{d_n}$, with $d_n=\dim(\mathbb{P}_n(\Omega)) = (n+1)(n+2)(n+3)/6$. Then the well-known orthogonal basis $\{\Phi_{\ell}\}_{\ell=1}^{d_n}$, with respect to the weight function $w({\mathbf{x}})$ on $[-1,1]^3$, consists of the tensor product of Chebyshev polynomials with total degree at most $n$, that is
\begin{eqnarray*}
    {\tilde{T}}_{\ell_1}(x_1) {\tilde{T}}_{\ell_2}(x_2)  {\tilde{T}}_{\ell_3}(x_3), {\hspace{0.2cm}} \ell=(\ell_i)_{i=1,2,3}, {\hspace{0.2cm}} \ell_1+\ell_2+\ell_3 \leq n,
\end{eqnarray*}
where ${\tilde{T}}_{k}(t)=\sqrt{2} \cos(k\arccos(t))$, for $k > 0$, ${\tilde{T}}_{0}(t) = 1$, $t \in [-1,1]$. 

For the quadrature rule on the cube, we use the method introduced in {\cite{DeMarchi2009new}}, which is based on a set of quadrature nodes $\mathcal{X}_{N_{2n}}$ and a weight function with the algebraic degree of exactness $2n$,  where the cardinality is
$N_{2n}  \approx \frac{(n+2)^3}{4}$. We claim that the hyperinterpolant will not be in general interpolant in the point set determined by the nodes because of $d_n < N_{2n}$.

This formula is determined as follows. Let $C_n=\{\cos(k\pi/n)\}_{k=0}^{n}$ be the set of $n+1$ Chebyshev-Lobatto points and let $C^E_{n+1}$, $C^O_{n+1}$ be, respectively, the restriction of $C_{n+1}$ to even and odd indices. Indeed, for any quadrature node ${\mathbf{\xi}}$, the corresponding weight $w_{\mathbf{\xi}}$ is
\begin{eqnarray*}
    w_{\mathbf{\xi}}:={\frac{4}{(n+1)^3}}\cdot \begin{cases}
        1, & {\mbox{if }} {\mathbf{\xi}} {\mbox{ is an interior point}}, \\
1/2, & {\mbox{if }} {\mathbf{\xi}} {\mbox{ is a face point}}, \\
1/4, & {\mbox{if }} {\mathbf{\xi}} {\mbox{ is an edge point}}, \\
1/8, & {\mbox{if }} {\mathbf{\xi}} {\mbox{ is a vertex point}}.
    \end{cases}
\end{eqnarray*}
Consequently, if $f,g \in {\mathbb{P}}_n(\Omega)$, then $\langle f,g \rangle=\langle f,g \rangle_{N_{2n}}$, where $\langle f,g \rangle_{N_{2n}}$ is the discrete scalar product defined by the such quadrature rule with exactness $2n$.

\begin{figure}[htbp]
  \centering
  \includegraphics[width=\textwidth]{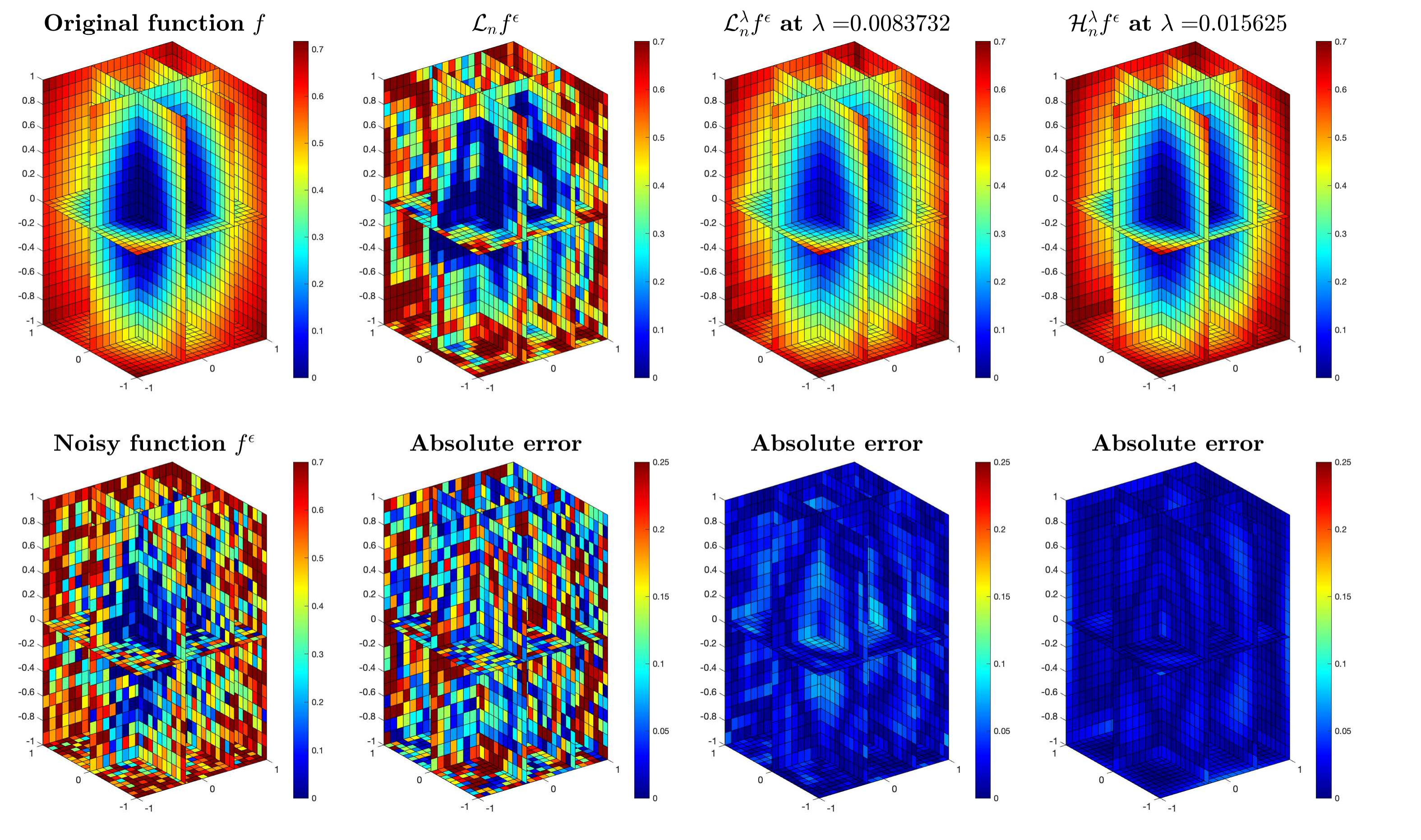}\\
  \caption{Approximate $f(x,y,z)=\text{exp}(-1/(x^2+y^2+z^2))$, perturbed by Gaussian noise ($\sigma = 0.2$), over $[-1,1]^3$, via  hyperinterpolation $\mathcal{L}_n f^{\epsilon}$,
    Lasso hyperinterpolation $\mathcal{L}_{n}^{\lambda}f^{\epsilon}$, and
    hard thresholding hyperinterpolation $\mathcal{H}_{n}^{\lambda} f^{\epsilon}$ at $n=20$.}\label{Cube_denoise}
\end{figure}

As observed in {\cite{DeMarchi2009new}}, setting
\begin{eqnarray*}
F({\mathbf{\xi}})=F(\xi_1,\xi_2,\xi_3)=\begin{cases}
        w_{\mathbf{\xi}} f({\mathbf{\xi}}), & {\mathbf{\xi}} \in \mathcal{X}_{N_{2n}}, \\
0, & {\mathbf{\xi}} \in (C_{n+1} \times C_{n+1} \times C_{n+1} ) \backslash {\hspace{0.05cm}} \mathcal{X}_{N_{2n}},
    \end{cases}
\end{eqnarray*}
the hyperinterpolation coefficients $\alpha_{\ell}$ are
\begin{eqnarray*}
    \alpha_{\ell} =\gamma_{\ell} \sum_{i=0}^{n+1} \left( \sum_{j=0}^{n+1} \left( \sum_{k=0}^{n+1} F_{ijk} \cos\frac{k {{\ell}_1} \pi}{n+1} \right) \cos \frac {j {{\ell}_2} \pi}{n+1}  \right) \cos \frac {i {{\ell}_3} \pi}{n+1},
\end{eqnarray*}
where $F_{ijk}=F\left( \cos\frac{i \pi}{n+1} , \cos\frac{j \pi}{n+1} ,\cos\frac{k \pi}{n+1}   \right)$, $i,j,k \in \{0,1,\ldots,n+1\}$ and
\begin{eqnarray*}
    \gamma_{\ell}=\prod_{s=1}^3 \gamma_{\ell_s},  {\hspace{0.2cm}}  \gamma_{\ell_s}=\begin{cases}
        \sqrt{2}, & \ell_s > 0,\\
1, & \ell_s = 0,
    \end{cases}
    {\hspace{0.2cm}} s=1,2,3.
\end{eqnarray*}
In view of this peculiar structure, fast computation of hyperinterpolation coefficients is feasible via FFT \cite{DeMarchi2009new}.

In our numerical examples, we examine the case of the function
\begin{eqnarray*}
    f(x,y,z)=\exp(-1/(x^2+y^2+z^2))
\end{eqnarray*}
contaminated by Gaussian noise $\mathcal{N}(0,\sigma^2)$ ($\sigma=0.2$).

\begin{figure}[htbp]
  \centering
  \includegraphics[width=\textwidth]{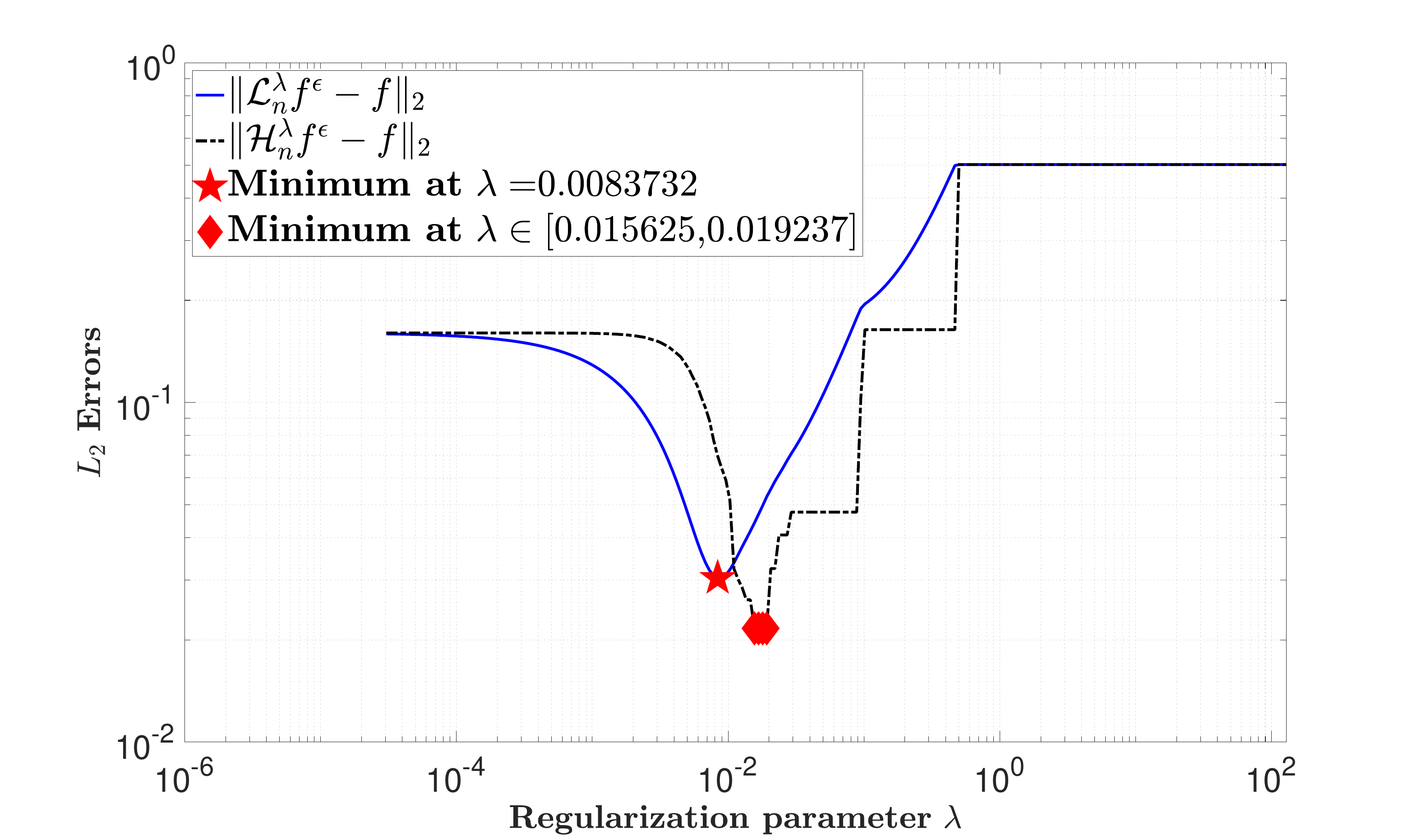}\\
  \caption{The choices of regularization parameter $\lambda$  for Lasso hyperinterpolation $\mathcal{L}_n^{\lambda} f^{\epsilon}$ and hard thresholding hyperinterpolation   $\mathcal{H}_n^{\lambda} f^{\epsilon}$ at $n=20$ approximating $f(x,y,z)=\text{exp}(-1/(x^2+y^2+z^2))$, perturbed by Gaussian noise ($\sigma = 0.2$), over $[-1,1]^3$.}\label{Cube}
\end{figure}

In Figures \ref{Cube_denoise}, we continue our investigation with $n=20$ and explore the denoising capabilities of the three hyperinterpolation methods. The findings align with the previous sections \ref{sec:sphere} and \ref{sec:spt}, where $\mathcal{H}_n^{\lambda} f^{\epsilon}$ demonstrates superior performance, as illustrated in Figure \ref{Cube}. Here, $\mathcal{H}_n^{\lambda} f^{\epsilon}$ achieves an $L_2$ error of 0.021623 at $\lambda \in [0.015625,0.019237]$ with 7 non-zero coefficients, while $\mathcal{L}_n^{\lambda} f^{\epsilon}$ reaches a minimum $L_2$ error of 0.030396 at $\lambda=0.0083732$ with 53 non-zero coefficients.
%In Figures \ref{Cube_denoise}  and \ref{Cube}, we fix the quadrature exactness $2n=40$ and thus have $N=2662$ cubature points. The denoising abilities of hard thresholding hyperinterpolation $\mathcal{H}_n^{\lambda}f^{\epsilon}$, Lasso hyperinterpolation $\mathcal{L}_n^{\lambda}f^{\epsilon}$ with regularization parameter $\lambda$ to reach their minimum $L_2$ errors, and classical hyperinterpolation $\mathcal{L}_n f^{\epsilon}$ are shown in Figure \ref{Cube_denoise}.   In Figure \ref{Cube}, $\mathcal{H}_n^{\lambda}f^{\epsilon}$ arrives at the lowest $L_2$ error of 0.021623 at $\lambda \in [0.015625,0.019237]$ with only 7 non-zero coefficients, while $\mathcal{L}_n^{\lambda} f^{\epsilon}$ reaches a minimum $L_2$ errors of 0.030396 at $\lambda=0.0083732$ but with 53 non-zero coefficients.

\section{Final remarks}
Concisely speaking, hard thresholding hyperinterpolation is the unique solution to an $\ell_0$-regularized weighted discrete least square problem, and has been proved to be an effective tool in denoising from the numerical examples. Hard thresholding hyperinterpolation satisfies the Pythagorean theorem with respect to the discrete (semi) inner product, which is an important geometric property that Lasso hyperinterpolation \cite{an2021lasso} and hybrid hyperinterpolation \cite{an2023hybrid} do not possess. Then we use the reciprocal of Christoffel function to prove that the upper bound of the uniform norm of hard thresholding hyperinterpolation operator is not greater than that of hyperinterpolation operator. What's more, a practical criterion, using the sum of the difference between the regularization parameter and the product of noise coefficients and signs of hyperinterpolation coefficients, is established to judge the denoising abilities of hard thresholding hyperinterpolation and Lasso hyperinterpolation. 

With the aid of the Marcinkiewicz-Zygmund property \cite{an2022exactness}, one can bypass the quadrature exactness as in \cite{an2024bypassing}, which can break the restriction of the application of hard thresholding hyperinterpolation. Once the quadrature exactness is not required, there are many quadrature rules that we can take, such as (Quasi) Monte-Carlo rules \cite{DickKuoSloan2013QMC}.  Furthermore, one may combine the springback penalty \cite{an2022springback} with the weighted least squares problem \eqref{equ:approximationproblem} to obtain a more stable and effective approximation scheme. In addition, it seems promising to discuss the relation between different types of noise and denoising ability of hard thresholding hyperinterpolation. 

\section*{Acknowledgements}
The first author (C. An) of the research is partially supported by National Natural Science Foundation of China (No. 12371099) and Tianfu Emei talent plan (No.1914).  The second author (J. Ran) gratefully acknowledges financial support from China Scholarship Council.

We are deeply grateful to the anonymous referees for their careful reading of our manuscript and their many insightful comments. We express our sincere thanks to Prof. Alvise Sommariva of the University of Padova for his helpful suggestions. We would also like to thank Dr. Hao-Ning Wu of the University of Georgia and Prof. Xiaoming Yuan of the University of Hong Kong for their continued interest in this work.

\section*{Conflict of interest}
The authors declare that they have no conflict of interest.

%\section*{References}
\bibliographystyle{siamplain}
\bibliography{HardHyperRef_JSC}

\begin{thebibliography}{10}

\bibitem{an2023hybrid}
{\sc C.~An, J.~Ran, and A.~Sommariva}, {\em {Hybrid hyperinterpolation over
  general regions}}, arXiv,  (2024),
  \href{http://dx.doi.org/https://doi.org/10.48550/arXiv.2305.05863}{doi:\nolinkurl{https://doi.org/10.48550/arXiv.2305.05863}}.

\bibitem{an2021lasso}
{\sc C.~An and H.-N. Wu}, {\em {Lasso hyperinterpolation over general
  regions}}, SIAM J. Sci. Comput., 43 (2021), pp.~A3967--A3991,
  \href{http://dx.doi.org/https://doi.org/10.1137/20M137793X}{doi:\nolinkurl{https://doi.org/10.1137/20M137793X}}.

\bibitem{an2022exactness}
{\sc C.~An and H.-N. Wu}, {\em {On the quadrature exactness in
  hyperinterpolation}}, BIT Numer. Math.,  (2022),
  \href{http://dx.doi.org/https://doi.org/10.1007/s10543-022-00935-x}{doi:\nolinkurl{https://doi.org/10.1007/s10543-022-00935-x}}.

\bibitem{an2024bypassing}
{\sc C.~An and H.-N. Wu}, {\em {Bypassing the quadrature exactness assumption
  of hyperinterpolation on the sphere}}, J. Complex., 80 (2024), p.~101789,
  \href{http://dx.doi.org/https://doi.org/10.1016/j.jco.2023.101789}{doi:\nolinkurl{https://doi.org/10.1016/j.jco.2023.101789}}.

\bibitem{an2023efficient}
{\sc C.~An and H.-N. Wu}, {\em {Is hyperinterpolation efficient in the
  approximation of singular and oscillatory functions?}}, J. Approx. Theory,
  299 (2024), p.~106013,
  \href{http://dx.doi.org/https://doi.org/10.1016/j.jat.2023.106013}{doi:\nolinkurl{https://doi.org/10.1016/j.jat.2023.106013}}.

\bibitem{an2022springback}
{\sc C.~An, H.-N. Wu, and X.~Yuan}, {\em The springback penalty for robust
  signal recovery}, Appl. Comput. Harmon. Anal., 61 (2022), pp.~319--346,
  \href{http://dx.doi.org/https://doi.org/10.1016/j.acha.2022.07.002}{doi:\nolinkurl{https://doi.org/10.1016/j.acha.2022.07.002}}.

\bibitem{atkinson2012spherical}
{\sc K.~Atkinson and W.~Han}, {\em {Spherical Harmonics and Approximations on
  the Unit Sphere: An Introduction}}, Lecture Notes in Mathematics, Springer,
  2012,
  \href{http://dx.doi.org/https://doi.org/10.1007/978-3-642-25983-8}{doi:\nolinkurl{https://doi.org/10.1007/978-3-642-25983-8}}.

\bibitem{bruckstein2009sparse}
{\sc A.~M. Bruckstein, D.~L. Donoho, and M.~Elad}, {\em From sparse solutions
  of systems of equations to sparse modeling of signals and images}, SIAM
  Review, 51 (2009), pp.~34--81,
  \href{http://dx.doi.org/https://doi.org/10.1137/060657704}{doi:\nolinkurl{https://doi.org/10.1137/060657704}}.

\bibitem{caliari2007hyperinterpolation}
{\sc M.~Caliari, S.~De~Marchi, and M.~Vianello}, {\em {Hyperinterpolation on
  the square}}, J. Comput. Appl. Math., 210 (2007), pp.~78--83,
  \href{http://dx.doi.org/https://doi.org/10.1016/j.cam.2006.10.058}{doi:\nolinkurl{https://doi.org/10.1016/j.cam.2006.10.058}}.

\bibitem{caliari2008hyperinterpolation}
{\sc M.~Caliari, S.~De~Marchi, and M.~Vianello}, {\em {Hyperinterpolation in
  the cube}}, Comput. Math. Appl., 55 (2008), pp.~2490--2497,
  \href{http://dx.doi.org/https://doi.org/10.1016/j.camwa.2007.10.003}{doi:\nolinkurl{https://doi.org/10.1016/j.camwa.2007.10.003}}.

\bibitem{dai2006hyperinterpolation}
{\sc F.~Dai}, {\em {On generalized hyperinterpolation on the sphere}}, Proc.
  Amer. Math. Soc., 134 (2006), pp.~2931--2941,
  \href{http://dx.doi.org/https://doi.org/10.1090/S0002-9939-06-08421-8}{doi:\nolinkurl{https://doi.org/10.1090/S0002-9939-06-08421-8}}.

\bibitem{sommariva2014multivariate}
{\sc S.~De~Marchi, A.~Sommariva, and M.~Vianello}, {\em {Multivariate
  Christoffel functions and hyperinterpolation}}, Dolomites Res. Notes Approx.,
  7 (2014), pp.~26--33,
  \href{http://dx.doi.org/https://doi.org/10.14658/PUPJ-DRNA-2014-Special_Issue-6}{doi:\nolinkurl{https://doi.org/10.14658/PUPJ-DRNA-2014-Special_Issue-6}}.

\bibitem{DeMarchi2009new}
{\sc S.~De~Marchi, M.~Vianello, and Y.~Xu}, {\em New cubature formulae and
  hyperinterpolation in three variables}, BIT Numer. Math., 49 (2009),
  pp.~55--73,
  \href{http://dx.doi.org/https://doi.org/10.1007/s10543-009-0210-7}{doi:\nolinkurl{https://doi.org/10.1007/s10543-009-0210-7}}.

\bibitem{delsarte1977spherical}
{\sc P.~Delsarte, J.~M. Goethals, and J.~J. Seidel}, {\em Spherical codes and
  designs}, Geom. Dedicata, 6 (1977), pp.~363--388,
  \href{http://dx.doi.org/https://doi.org/10.1007/BF03187604}{doi:\nolinkurl{https://doi.org/10.1007/BF03187604}}.

\bibitem{DickKuoSloan2013QMC}
{\sc J.~Dick, F.~Y. Kuo, and I.~H. Sloan}, {\em {High-dimensional integration:
  The quasi-Monte Carlo way}}, Acta Numer., 22 (2013), pp.~133--288,
  \href{http://dx.doi.org/https://doi.org/10.1017/S0962492913000044}{doi:\nolinkurl{https://doi.org/10.1017/S0962492913000044}}.

\bibitem{donoho1994ideal}
{\sc D.~L. Donoho and I.~M. Johnstone}, {\em {Ideal spatial adaptation by
  wavelet shrinkage}}, Biometrika, 81 (1994), pp.~425--455,
  \href{http://dx.doi.org/https://doi.org/10.1093/biomet/81.3.425}{doi:\nolinkurl{https://doi.org/10.1093/biomet/81.3.425}}.

\bibitem{Foucart2013compressing}
{\sc S.~Foucart and H.~Rauhut}, {\em A Mathematical Introduction to Compressive
  Sensing}, Applied and Numerical Harmonic Analysis, Birkh{\"a}user New York,
  NY, 1~ed., 2013,
  \href{http://dx.doi.org/https://doi.org/10.1007/978-0-8176-4948-7}{doi:\nolinkurl{https://doi.org/10.1007/978-0-8176-4948-7}}.

\bibitem{hansen2009norm}
{\sc O.~Hansen, K.~Atkinson, and D.~Chien}, {\em {On the norm of the
  hyperinterpolation operator on the unit disc and its use for the solution of
  the nonlinear Poisson equation}}, IMA J. Numer. Anal., 29 (2009),
  pp.~257--283,
  \href{http://dx.doi.org/https://doi.org/10.1093/imanum/drm052}{doi:\nolinkurl{https://doi.org/10.1093/imanum/drm052}}.

\bibitem{MR2274179}
{\sc K.~Hesse and I.~H. Sloan}, {\em Hyperinterpolation on the sphere}, in
  Frontiers in interpolation and approximation, vol.~282 of Pure Appl. Math.,
  Chapman \& Hall/CRC, Boca Raton,, 2007, pp.~213--248.

\bibitem{LeGia2001uniform}
{\sc Q.~T. Le~Gia and I.~H. Sloan}, {\em The uniform norm of hyperinterpolation
  on the unit sphere in an arbitrary number of dimensions}, Constr. Approx., 17
  (2001), pp.~249--265,
  \href{http://dx.doi.org/https://doi.org/10.1007/s003650010025}{doi:\nolinkurl{https://doi.org/10.1007/s003650010025}}.

\bibitem{lin2021distributed}
{\sc S.-B. Lin, Y.~G. Wang, and D.-X. Zhou}, {\em {Distributed filtered
  hyperinterpolation for noisy data on the sphere}}, SIAM J. Numer. Anal., 59
  (2021), pp.~634--659,
  \href{http://dx.doi.org/https://doi.org/10.1137/19M1281095}{doi:\nolinkurl{https://doi.org/10.1137/19M1281095}}.

\bibitem{Nevai1986Christoffel}
{\sc P.~Nevai}, {\em {G{\'e}za Freud, orthogonal polynomials and Christoffel
  functions. A case study}}, J. Approx. Theory, 48 (1986), pp.~3--167,
  \href{http://dx.doi.org/https://doi.org/10.1016/0021-9045(86)90016-X}{doi:\nolinkurl{https://doi.org/10.1016/0021-9045(86)90016-X}}.

\bibitem{reimer2003multivariate}
{\sc M.~Reimer}, {\em {Multivariate Polynomial Approximation}}, Birkh{\"a}user
  Basel, first~ed., 2003,
  \href{http://dx.doi.org/https://doi.org/10.1007/978-3-0348-8095-4}{doi:\nolinkurl{https://doi.org/10.1007/978-3-0348-8095-4}}.

\bibitem{sloan1995hyperinterpolation}
{\sc I.~H. Sloan}, {\em {Polynomial interpolation and hyperinterpolation over
  general regions}}, J. Approx. Theory, 83 (1995), pp.~238--254,
  \href{http://dx.doi.org/https://doi.org/10.1006/jath.1995.1119}{doi:\nolinkurl{https://doi.org/10.1006/jath.1995.1119}}.

\bibitem{sloan1997interpolation}
{\sc I.~H. Sloan}, {\em Interpolation and hyperinterpolation on the sphere}, in
  Multivariate Approximation: Recent Trends and Results, W.~Haussmann,
  K.~Jetter, and M.~Reimer, eds., vol.~101 of Mathematical Research, Akademie
  Verlag GmbH, Berlin (Wiley-VCH), 1997, pp.~255--268.

\bibitem{sloan2012filtered}
{\sc I.~H. Sloan and R.~S. Womersley}, {\em {Filtered hyperinterpolation: a
  constructive polynomial approximation on the sphere}}, Int. J. Geomath., 3
  (2012), pp.~95--117,
  \href{http://dx.doi.org/https://doi.org/10.1007/s13137-011-0029-7}{doi:\nolinkurl{https://doi.org/10.1007/s13137-011-0029-7}}.

\bibitem{Sommariva2021Tchakaloff}
{\sc A.~Sommariva and M.~Vianello}, {\em {Near-algebraic Tchakaloff-like
  quadrature on spherical triangles}}, Appl. Math. Lett., 120 (2021),
  p.~107282,
  \href{http://dx.doi.org/https://doi.org/10.1016/j.aml.2021.107282}{doi:\nolinkurl{https://doi.org/10.1016/j.aml.2021.107282}}.

\bibitem{sommariva2021spheritri}
{\sc A.~Sommariva and M.~Vianello}, {\em Numerical hyperinterpolation over
  spherical triangles}, Math. Comput. Simul., 190 (2021), pp.~15--22,
  \href{http://dx.doi.org/https://doi.org/10.1016/j.matcom.2021.05.003}{doi:\nolinkurl{https://doi.org/10.1016/j.matcom.2021.05.003}}.

\bibitem{stein2011fourier}
{\sc E.~M. Stein and R.~Shakarchi}, {\em {Fourier analysis: an introduction}},
  vol.~1, Princeton University Press, 2011.

\bibitem{Wade2013hyperinterpolation}
{\sc J.~Wade}, {\em On hyperinterpolation on the unit ball}, J. Math. Anal.
  Appl., 401 (2013), pp.~140--145,
  \href{http://dx.doi.org/https://doi.org/10.1016/j.jmaa.2012.11.052}{doi:\nolinkurl{https://doi.org/10.1016/j.jmaa.2012.11.052}}.

\bibitem{Wang2017needlet}
{\sc Y.~G. Wang, Q.~T. Le~Gia, I.~H. Sloan, and R.~S. Womersley}, {\em Fully
  discrete needlet approximation on the sphere}, Appl. Comput. Harmon. Anal.,
  43 (2017), pp.~292--316,
  \href{http://dx.doi.org/https://doi.org/10.1016/j.acha.2016.01.003}{doi:\nolinkurl{https://doi.org/10.1016/j.acha.2016.01.003}}.

\bibitem{womersley2018efficient}
{\sc R.~S. Womersley}, {\em Efficient {S}pherical {D}esigns with {G}ood
  {G}eometric {P}roperties}, in Contemporary {C}omputational {M}athematics -
  {A} {C}elebration of the 80th {B}irthday of {I}an {S}loan, Springer
  International Publishing, Cham, 2018, pp.~1243--1285,
  \href{http://dx.doi.org/https://doi.org/10.1007/978-3-319-72456-0_57}{doi:\nolinkurl{https://doi.org/10.1007/978-3-319-72456-0_57}}.

\end{thebibliography}
\clearpage

\end{document}